\documentclass[10pt]{article}

\usepackage{amsmath,latexsym,amssymb,amsthm}

\numberwithin{equation}{section} \allowdisplaybreaks

\newtheorem{theorem}{\bf\normalsize Theorem}[section]
\newtheorem{corollary}[theorem]{\bf\normalsize Corollary}

\newtheorem{lemma}[theorem]{\bf\normalsize Lemma}

\newtheorem{remark}[theorem]{\bf\normalsize Remark}

\begin{document}

\title{A note on    Bourgain's slicing problem}
\author{QingYang Guan\footnote{
 Institute of Applied Mathematics,  AMSS, CAS.
   }
   }
\date{} \maketitle

\begin{abstract}
This note is to study    Bourgain's slicing problem   following the    routes investigated in the last decade. We show that the
      slicing
constant $L_n$ is bounded by $C\log(\log n) $, $n\geq 3$, for some universal constant $C$.
\end{abstract}


{\textbf{\small\noindent AMS(2000) Subject Classification:
  Primary $52\mathrm{A}23$;\  Secondary
  $60\mathrm{H}30$ }

\section{Introduction}
\subsection{\normalsize the main results}
The main subject of this note  is to give a  iterated logarithm  bound for   the  constant in  Bourgain's slicing problem.       The Bourgain's slicing problem is one of the central questions in high-dimensional convex geometry. It is about  a universal regularity  of the  convex sets in
$\mathbb{R}^n$. A standard form of this problem is   whether  for any convex set $K\subseteq \mathbb{R}^n$ with volume one, there exists a hyperplane $H\subseteq \mathbb{R}^n$ such that
$$
\mathrm{Vol}_{n-1}(K\cap H)\geq c,\  \ \ \ n\geq 2,
$$
where $c>0$ is a universal constant, and $ \mathrm{Vol}_{n-1}$
refers to  the standard  $(n-1)$-dimensional volume measure on $H$.
This motivates us to estimate  the constant      $L_n$ defined by
$$
(L_n)^{-1}=  \inf_{K\subseteq \mathbb{R}^n} \sup_{H\subseteq \mathbb{R}^n}\mathrm{Vol}_{n-1}(K\cap H),\  \ \ \ n\geq 2,
$$
where $H$  runs over all hyperplanes and  $K$ runs over all convex sets with volume one. The main result of this note is the  following   bound of $L_n$.

 \begin{theorem}\label{theorem}
There exists a universal constant $C>0$ such that $L_n \leq C\log(\log n)$
for   $n\geq 3$.
\end{theorem}

Applying a well known connection   established  in
  Eldan and Klartag \cite{EK11}, we prove the bound above  via proving the same kind  of   bound for the   thin cell constant $\sigma_n$.
The thin cell parameter  $\sigma_\mu$ of a log-concave probability measure $\mu$ is defined by
  \begin{align}
 n\sigma_\mu^2= \mathrm{Var}_\mu(|x|^2)  ,
 \end{align}
 where $\mathrm{Var}_\mu (u)=\int u^2 d\mu-(\int u d\mu)^2$ is the variance of a function $u$.
 Set
  \begin{align}
\sigma_n= \sup \{\sigma_\mu:\ \mu\ \mathrm{is\ an\ isotropic\    \textrm{log-concave} \ probability\ measure\ on}\   \mathbb{R}^n \}  .\nonumber
 \end{align}
 Here $\mu$  is called isotropic
if its barycenter is zero  and its  covariance matrix is the identity.
The   connection between  $L_n$ and $\sigma_n$  in
\cite{EK11}  is
 \begin{align}\label{con}
L_n \leq C\sigma_n ,\ \ \ \ n\geq 1,
 \end{align}
 where $C>0$ is a universal  constant.

 Whether $\sigma_n$ can be bounded by a universal constant is called the  thin cell conjecture. The thin cell conjecture
comes from  the works of  Anttila, Ball and Perissinaki
\cite{ABP03} and Bobkov and Koldobsky \cite{BK03}, which is   later than
the   slicing problem  originated  in Bourgain
\cite{B86,B87}. There is still a more stronger  conjecture called the KLS
isoperimetric conjecture  given in   Kannan,
Lov$\acute{\mathrm{a}}$sz, and Simonovits \cite{KLS95}.
 More and more connections have been discovered gradually among these three conjectures although they may  have  different    motivations.

The isoperimetric constant of $\mu$ with density $\rho$ is define via
 \begin{align}
   \frac{1}{\psi_\mu} =\inf_{A\subseteq \mathbb{R}^n}\Big\{\frac{\int_{\partial A} \rho }{min\{\mu(A),1-\mu(A)\}}\Big\} ,\nonumber
   \end{align}
where $A$ runs over all open sets with smooth boundary satisfying  $0<\mu(A)<1$.
Set $\psi_n= \sup_{\mu} \psi_\mu$, where $\mu$ runs over all isotropic, log-concave probability measure on $\mathbb{R}^n$.
The   KLS
isoperimetric conjecture states  that $\psi_n\leq C$ for a   universal constant  $C>0$.

The best bound of $L_n$
 before this note  is $C\sqrt{\log n}$, $n\geq 2$, which is a consequence of (\ref{con}) and  the recent breakthrough    $\psi_n\leq C\sqrt{\log n}$ in   Klartag \cite{K23}.
The bound of $L_n$ given in  Theorem  \ref{theorem} above  is a direct consequence of   this isoperimetric  constant  estimate, (\ref{con}) and the following relation.
 \begin{theorem}\label{theorem2}
There exists a universal constant $C>0$ such that $\sigma_n \leq C\log  \psi_n$.
\end{theorem}
Combing (\ref{con}) and Theorem 1.2, we get  $L_n \leq C\log  \psi_n$, which  improves the previous  results of $L_n \leq    e^{C\psi_n^2}$ in Ball and   Nguyen \cite{BN12} and $L_n \leq     C\psi_n$
by Eldan and Klartag \cite{EK11}. This result has an unexpected character  that
a weak control of the isoperimetric  constant  $\psi_n$ by the slicing constant $L_n$ or the thin cell constant $\sigma_n$ would lead to the truth of these three  conjectures.

\subsection{\normalsize a short introduction of the previous  development}

Next we briefly introduce some  basic concepts,     methods and
some    development of the subject investigated  before this note. Definitions and notations below will be used in the main content.
Although  we mainly focus  on the slicing conjecture,   the statement  has closed connection with the other two conjectures mentioned  above.
 We refer
to Klartag and Milman \cite{KM22} and the references therein for
further backgrounds,  equivalent formulations and applications of the
problem.
For a more complete  statement,
we refer to  Eldan and Klartag  \cite{EK11},  Eldan \cite{E13}, Lee and  Vempala \cite{LV17}\cite{LV24}, Chen \cite{C21}, Klartag and  Putterman \cite{KP21}, Klartag   Lehec \cite{KL22},
Jambulapati,  Lee and  Vempala \cite{JLV22},
Klartag \cite{K21}\cite{K23} and the references therein. See also Klartag and  Lehec \cite{KL24} for a   general introduction beyond the subject of  this note.
Notations here are mainly taken from \cite{K23}, \cite{KP21} and \cite{KL22}.

$\mathbf{A})$ \ Let $\mu$ be a log-concave probability measure with    density  $\rho$.
   In the last ten years, a driven force in the  study of the  high-dimensional convex sets  is  the   Eldan's
stochastic localization started in  Eldan \cite{E13},  which is  a kind of   partition of   $\mu$  through a random procedure.
     Unlike the
determined localization, in which the dimension of a
fragment  finally goes to  one,   this stochastic version
  is   a kind of partition
from dimension $n$ to dimension $n$, in  which a fragment  may keep more information of the original one than the previous method.  In fact, a fragment in
this stochastic    partition is the original measure $\mu$ weighted by
a Gaussian function, which makes it even more regular than $\mu$ in some aspect.

More concretely, the   Eldan's stochastic localization    of a probability measure $\mu$  is a probability distribution  valued stochastic process
$(p_t)_{t\geq 0}$.
Here    $p_t$ refers to the density of the corresponding distribution and  $p_0=\rho$.
  In
a simplified  version by Lee and Vemplala \cite{LV17}, it is the solution of  the following stochastic differential  equation with the initial condition $p_0$:
\begin{align}
  dp_t(x)= p_t(x) \langle x-a_t,d B_t\rangle ,\ \ \ \ \ \   x\in \mathbb{R}^n,\label{p}
 \end{align}
where   $a_t $ is the barycenter of $p_t$ and $(B_t)=(B_{k,t})_{1\leq k\leq n}$ is  a standard Brownian motion in $\mathbb{R}^n$.
Denote, respectively,  the components of  the barycenter $a_t $ and the covariance matrix $A_t$ of $p_t$ by
\begin{align}
  a_{i;t}=& \int
    x_i p_t(x)dx ,\ \ \ \ \ \ \ \ \ \ \ \ \ \ \ \ \ \ \ \ \ \ \ \ \ \ \ 1\leq i\leq n,\nonumber\nonumber\\
  A_{ij;t}=&\int  (x_i-a_{i;t})(x_j-a_{j;t})p_t(x) dx,\ \ \ \ 1\leq i,j\leq n.\nonumber
 \end{align}
 In \cite{LV17}, it follows from  (\ref{p}) that
\begin{align}
 d a_t= \int  x \langle x-a_t,d B_t\rangle p_t(x)dx=A_tdB_t .\nonumber
 \end{align}
In Klartag and  Lehec \cite{KL22}, it shows further by  It$\mathrm{\hat{o}}$ formula  that
\begin{align}\label{one dir}
 \frac{d }{dt}\mathbb{E}|a_t|^2=  \mathbb{E}\mathrm{Tr}( A_t^2),\ \ \ \ t>0 .
 \end{align}

As the original version in    Eldan \cite{E13}, there is another    way to construct   $(p_t)$, which involves the Gaussian function mentioned above. For $t\geq 0$ and
$\theta\in \mathbb{R}^n$, define the
probability density function
\begin{align}\label{pt}
p_{t,\theta}(x)=Z(t,\theta)^{-1}e^{\langle
\theta,x\rangle-\frac{t|x|^2}{2}}p_0(x),\ \ \ \ \ \ x\in \mathbb{R}^n,
\end{align}
where $Z(t,\theta)=\int e^{\langle
\theta,x\rangle-\frac{t|x|^2}{2}}p_0(x)dx$.
Set   $a(t,\theta)=\int x  p_{t,\theta}(x)dx$. In   Lee and Vemplala \cite{LV17}, it shows that
 $(p_{t,\theta_t})$ is a solution of (\ref{p}) when   $(\theta_t)$ is  the solution of the   stochastic differential equation
\begin{align}
  d{\theta}_t= a(t,\theta_t)dt+dB_t,\ \ \ \theta_0=0.\nonumber
\end{align}

In Klartag and  Putterman \cite{KP21}, it shows  further that  the stochastic process
$(\theta_t)$ is the same as the   process
$( tX+B_t )$ in the sense of distribution. Here $X$ is  a random vector with distribution $\mu$,  independent of $(B_t)$.
Moreover,  $p_{t,\theta}$ can be taken as the  density    of $X$ under the condition that  $tX+B_t=\theta$.

$\mathbf{B}) $\  Notice that $(p_t)$ does not have a semigroup structure on a finite dimensional space. In Klartag and  Putterman \cite{KP21}, the semigroup picture of $
 (p_{t,\theta})$ on $\mathbb{R}^n$ is revealed through  a transformation.
 Set $Q_0=Id$.  When  $u\in L^1(\mu)$, the  following relation
is given in \cite{KP21}:
    \begin{align}
  Q_s u(y)=\int u(x)p_{t,\theta}(x)dx,\ \ \ \ \ with\ s=1/t,\  y=\theta/t,\ \ s,t>0,\label{st}
 \end{align} where $(Q_s)$ represents part of  a nonhomogeneous heat semigroup which can   be constructed    independently from  $(p_{t,\theta})$.

Write $\mu_0=\mu$, $ \rho_0(x)=\rho(x)$ and  denote
  \begin{align}
 \mu_s=\mu \ast \gamma_s, \ \ \ \ \ d\mu_s=\rho_s(x)dx,\ \ \ \ s>0,\nonumber
 \end{align}
where $\gamma_s$ is the Gaussian measure on $\mathbb{R}^n$ with mean zero and covariance $s\cdot Id$.
 Let $(P_s)$ be the standard heat semigroup on $\mathbb{R}^n$, i.e., $P_0=\mathrm{Id}$ and
     \begin{align}
 P_s u=u\ast \gamma_s,\ \ \ \  s>0.\nonumber
 \end{align}
 In \cite{KP21}, it shows that   $P_s$  is a  contraction operator from $L^2(\mu_s)$ to $L^2(\mu)$, and $Q_s$ can be defined as  the adjoint operator of $P_s$  from   $L^2(\mu)$ to $L^2(\mu_s)$.  It has the following form
 \begin{align}
  Q_s u(y)= \frac{P_s(u\rho)}{\rho_s},\ \ \ s\geq 0.\label{Q}
 \end{align}
Since the distribution of $\theta_t/t$ is the same as that of $X+\frac{B_t}{t}$, for $u\in L^1(\mu)$ we   get by (\ref{st})
 \begin{align}
 \mathbb{E} \Big| \int udp_t\Big|^2=\int\big| Q_s u\big|^2 d\mu_s,\ \ \ \ \ \  where\ s=1/t,  \ s,t>0.\label{Qq}
 \end{align}
See Lemma 2.2 in  Klartag and  Lehec \cite{KL22} for more details.

 The differential of $(Q_s)$ with respect to  $s$ is the   box operator
  \begin{align}
\Box_s =\frac{\Delta}{2} +\nabla \log \rho_s \cdot\nabla.
   \nonumber
 \end{align}
 In  Klartag and  Putterman \cite{KP21}, with   operators $(\Box_s)$, the   $\Gamma$-calculus   of  Markov  semigroup (see e.g. \cite{BGL14}) is extended to this nonhomogeneous setting.
 For $s\geq 0$ and   functions $u,v$, which are smooth for example, set $\Gamma_{0,s}(u,v)=uv$ and define inductively
    \begin{align}
\Gamma_{k,s}(u,v)= \Box_s \Gamma_{k-1,s}(u,v)-
\Gamma_{k-1,s}(\Box_su,v)-\Gamma_{k-1,s}(u,\Box_sv)-\frac{d\Gamma_{k-1,s}}{ds}(u,v) ,\ \ \ k\geq1.
   \nonumber\end{align}
It shows in Proposition 2.3 \cite{KP21} that   for $ k= 0,1$,
 \begin{align}
 \frac{d}{ds}\int  \Gamma_{k,s}(Q_s u,Q_s v) d\mu_s  =-\int \Gamma_{k+1,s}(Q_s u,Q_s v) d\mu_s,\ \ s>0,
   \label{di}
 \end{align}
where $u$ and $v$ both have  subexponential decay  relative to $\rho$ of which the definition is     stated  in (\ref{se}) below.

$\mathbf{C})$ \ In this  part we turn to  the spectral analysis for the Laplacian  of $\mu$ which  is defined by
  \begin{align}
L  = \Delta +\langle \nabla \log \rho,\nabla\rangle  .\label{gen}
 \end{align}
  A crucial property of the formula (\ref{pt}) is that $p_t$ is a $t$-uniformly log-concave probability measure.
A   probability measure $\mu$ with density $\rho$ is called $t$-uniformly log-concave if   $\rho(x)=e^{-\varphi(x)-t|x|^2}$ for some $\mathbb{R}\cup \{+\infty\}$ valued convex function $\varphi$.
The log-concave probability measure is just the  $0$-uniformly log-concave  probability measure. It does not include the degenerate case here.

 Denote by $C_p(\mu)$ the Poincar$\acute{\mathrm{e}}$ constant of $\mu$ which is  the smallest constant $C$ such that for any locally Lipschitz function $u\in L^2(\mu)$,
  \begin{align}
 \mathrm{Var}_\mu(u)\leq  C\int  |\nabla u|^2 d\mu.\nonumber
 \end{align}
 When $\mu$ is  $t$-uniformly log-concave for some $t>0$, we have
the    log-concave Lichnerowicz inequality (see e.g. \cite{BGL14})
\begin{align}\label{Lic}
 C_p(\mu)\leq   \frac{1}{t} .
 \end{align}
This  estimate is also a consequence of  the Brascamp-Lieb inequality  in \cite{BL76}.
Applying     Bochner formula to an     eigenfunction corresponding to the   eigenvalue $C_p(\mu)^{-1}$,   (\ref{Lic}) is improved by the following
 log-concave spectral-variance  inequality
   in  Theorem  1.3.   Klartag   \cite{K23},
\begin{align}\label{kla}
 C_p(\mu)\leq \sqrt{ \frac{\|Cov(\mu)\|_{op} }{t}} ,\ \ \ \ \ t>0.
 \end{align}
  For an arbitrary  isotropic  log-concave probability measure, it may not be   $t$-uniformly log-concave for any $t>0$.
However,  when $\rho$ is positive and smooth,  in (\ref{theta}) below  we   show that it is $1$-uniformly log-concave in an average sense.

Another notion  related to the spectral of $L$ is the $H^{-1}(\mu)$ norm. For $u\in L^2(\mu)$ with $\int u d\mu=0$, define
 \begin{align}
  \| u\|_{H^{-1}(\mu)}=\sup\{\int uf d\mu:\ f\in L^2(\mu)\ is \ locally \ Lipshcitz\ with\ \int |\nabla f|^2 d\mu\leq 1 \}.\nonumber
   \end{align}
According to Proposition 10 in Barthe and  Klartag  \cite{BK20}, for any smooth function $u$ such that $u,|\nabla u|\in L^2(\mu)$ and $\int \nabla u d\mu=0$,
 \begin{align}
  Var_\mu(u)\leq  \sum_{i=1}^n\|\partial_i u\|_{H^{-1}(\mu)}^2 .\label{ct}
   \end{align}
The formula above is fundamental in Klartag \cite{K09} and Klartag and  Lehec \cite{KL22} to estimate  the thin cell parameter.

$\mathbf{D})$\
The bound of the slicing constant $L_n$ has been improved for several times. For example,
 $L_n\leq C n^{\frac{1}{4}}\log n$ in Bourgain \cite{B91}, $L_n\leq C n^{\frac{1}{4}} $ in Klartag \cite{K06},
$L_n\leq e^{C \sqrt{(\log n)( \log   \log n})} $ in Chen \cite{C21},  $L_n\leq C (\log n)^4$ in  Klartag and  Lehec \cite{KL22},
$L_n\leq C (\log n)^{2.223...}$ in   Jambulapati,  Lee and  Vempala \cite{JLV22} and $L_n\leq C\sqrt{ \log n}$ in Klartag \cite{K23}.
 For a historic comments on the thin  cell constant  $\sigma_n$, we refer to \cite{KL22}.

\subsection{\normalsize introduction of the proof}

As stated above, to prove Theorem \ref{theorem}, we only need to establish the relation in Theorem \ref{theorem2}.
In  Klartag and  Lehec \cite{KL22}, applying  the $H^{-1}$-inequality  (\ref{ct}) above, the estimate of the thin cell parameter $\sigma_\mu$ is reduced  to the estimate  of  the spectral measure  of the coordinate functions.
More explicitly, we need to control $M_\lambda:=\sum_{i=1}^n \int |E_\lambda x_i|^2 d\mu$, where  $E_\lambda$ is  the spectral projection of $-L$ below level $\lambda>0$.
In     \cite{KL22}, to estimate $M_\lambda$, a relation between $(M_\lambda)$ and $(N_s):=(\sum_{i=1}^n\|Q_s x_i\|_{L^2(\mu_s)}^2)$ is established  by spectral theory and $\Gamma$-calculus.
 Then, the rest part of the proof   is to give a suitable bound of $(N_s)$.
However, to  give  a iterated logarithm bound,  a stronger relation between $(M_\lambda)$ and $(N_s)$ is necessary.  Based on the $\Gamma$-calculus developed in Klartag and  Putterman \cite{KP21},
a new  connection between $(M_\lambda)$ and $(N_s)$ is formulated  in Lemma \ref{TecN} below.

It follows  from (\ref{Qq}) that   $N_s=\mathbb{E}|a_{1/s}|^2$. By (\ref{one dir}), to
estimate $\mathbb{E}|a_t|^2$  it is natural to  study   the
covariance matrix $A_t$, which  is a central object  in the study of  $(p_t)$.
The following is a partial  outline of the known results   of   $A_t$. The parameter $\kappa_n$ is stated   in (\ref{kappa}) below.  \medskip

 \emph{Small  time  $t\in (0, (C\kappa_n^2 \log n)^{-1} ]$}: \  \  $\mathbb{P}(\|A_t\|_{op}\geq 2)\leq e^{-(Ct)^{-1}}$   in Klartag and  Lehec \cite{KL22}; \medskip

\emph{global    time  $t> 0$}: \ \   \ \  \ \  \ \  \ \ \ \ \ \   \ \  \ \  \ \  \ \ $\|A_t\|_{op}\leq
\frac{1}{t}$ by (\ref{Lic});\medskip

\emph{global    time  $t> 0$}: \ \ \ \   \ \  \ \  \ \ \ \ \ \ \ \   \ \  \ \  \ \     for $q\geq 3$ and  $t>s>0$,
$\mathbb{E} \mathrm{Tr}(A_t^q)\leq
\big(\frac{t}{s}\big)^{2q}\mathbb{E}\mathrm{Tr}(A_s^q)$  in  Chen \cite{C21};\medskip

\emph{global    time  $t> 0$}: \ \  \ \   \ \  \ \  \ \  \ \ \ \  \ \  \ \  \ \  \ \   for $t>s>0$,
$\mathbb{E}\mathrm{Tr}(A_t^2)\leq
\big(\frac{t}{s}\big)^{3}\mathbb{E}\mathrm{Tr}(A_s^2)$ in  Klartag and  Lehec \cite{KL22};\medskip

\emph{global    time  $t> 0$}: \ \  \ \   \     for $t>s>0$,
$\mathbb{E}\mathrm{Tr}(A_t^2)\leq
\big(\frac{t}{s}\big)^{2\sqrt{2}}\mathbb{E}\mathrm{Tr}(A_s^2)$ in  Jambulapati,  Lee and  Vempala \cite{JLV22}.\medskip

Notice that the bound above  in Chen \cite{C21} is proved in  the Eldan's  original model. A similar result in the simplified  model by  Lee and Vemplala is given in Klartag \cite{K21}.
See also Eldan \cite{E13},   Lee and  Vempala\cite{LV17}\cite{LV24} for some earlier related results.

  In (\ref{000med}) below, we give a sharp  estimate  of $\mathbb{E}\mathrm{Tr}(A_t^2)$ up to a constant factor.
Instead of estimating  $\mathbb{E}\mathrm{Tr}(A_t^2)$ directly as before,  we consider some  assistant  stochastic processes in (\ref{ap}) stated briefly as below:
   \begin{align}
  F_{k;t}= \sum_{i=1}^nf_{k}(\lambda_{i;t}) ,\ \ \ \ \ t\geq 0,\   k\geq 1,\nonumber
   \end{align}
where $(\lambda_{i;t}) $ are    eigenvalues of the covariance matrix $A_t$.
Here, each of   $(f_k)$ is a positive and increasing function on  $ [0,\infty)$ with the  continuous second derivative. It is a square function  when the    independent variable
is relatively large
and a  exponential function when it is small. The purpose  of these functions is to restrict the
analysis of $A_t$ on the large  eigenvalues. As a consequence,  a bound of $\mathbb{E} \mathrm{Tr}(A_t^2)$ can be deduced from  an estimate of $\mathbb{E} F_{k;t}$.
 Meanwhile, it allows us to use the previous estimates for the operator
norm of $A_t$ at time $(\log n)^{-\alpha}$ for some $\alpha >0$. To achieve this aim, we need a finite sequence of $( F_{k;t})$ with length depending  on $n$.
A process  $ F_{k;t}$ is useful only  on a particular interval $[t_k,t_{k+1}]$. See  (\ref{intr}) for more details.

The construction of $(f_k)$ is  based on   Lemma \ref{Step5}. The estimate of $(\mathbb{E} F_{k;t})$ is based on  Lemma \ref{Step4}. An upper bound  and a lower
   bound  of $\mathbb{E} \mathrm{Tr}(A_t^2)$
deduced from
the estimates of $(\mathbb{E} F_{k;t})$  are  stated
in Lemma \ref{Step0}. The final  estimate of $\mathbb{E} |a_t|^2$ in the proof of Theorem \ref{theorem2} is optimal within a constant factor .

Throughout the paper, notations $n,i,j,k$ are always denoted for
    integers  with $n\geq 1$. Notation
   $x=(x_1,\cdots, x_n)$ is always denoted for
 elements  of $ \mathbb{R}^n$.
The letter $C$ denotes positive universal constant, whose value may change from one line to the next. Notation $\log$
stands for the natural logarithm.      For a $n$ by $n$ matrix $A$,   $\|A\|_{op}$ stands for  the   operator norm. We also denote the trace of $A$ by $\mathrm{Tr} (A)$.
For $n$ by $n$ symmetric  matrices  $A$ and $B$, write $A\geq B$ if  $A-B$ is positive semi-definite. Id is for the identity matrix. For a function $f$ defined on the real line, $f'$ and $f''$ stand for
its first derivative and the second derivative, respectively.  Notations $\nabla$, $\nabla^2$ refers to the gradient and Hessian operators, respectively.   Denote by $\langle \cdot,\cdot\rangle$ the
 standard inner product of $\mathbb{R}^n$. Notation $\mathbb{E}$ refers to the expectation operation and $\int$ is always for $\int_{\mathbb{R}^n}$
when the lower index is not specified. For a  set $H$,  the  indicator function of it is denoted by $I_H$.

\section{assistant function}

The aim of this section is to make some preparation to prove the following lemma in the next section.
\begin{lemma}\label{Step0}
 There exists a universal constant $C>0$ such that for any isotropic  log-concave probability measure $\mu$ on   $\mathbb{R}^n$ with $n>C$,
   \begin{align}\label{000med}
   \frac{n}{2}\leq \mathbb{E} \mathrm{Tr}(A_t^2)\leq 8n,\ \ \ \ \ t\in [0,e^{-10^3}] .
 \end{align}
\end{lemma}\medskip

 \begin{lemma}\label{Step5}
Let   $D_0>4$ and    $r_0\in[ \frac{7}{3},\frac{8}{3}]$.  Then there exists a constant $b$ and a
increasing function $f \in C^2(\mathbb{R})$ such that
   \begin{align}\label{fr}
   f(r)=&e^{D_0(r-r_0)},\ \ \ \ \ r\leq r_0-D_0^{-1};\\ \ \ \ \ f(r)=&br^2,\ \ \ \ \ \ \ \ \ \ \ \ r\geq r_0, \label{fr2}
   \end{align}
   and
     \begin{align}
  \frac{1}{20}\leq b\leq \frac{1}{5};\ \ \ \ \ \ \  |f''(r)|\leq   D_0^2f(r) ,\ \ \ \ r\in \mathbb{R}. \label{sece}
   \end{align}
\end{lemma}\noindent\textbf{Proof}\ Let $D_0>4$ and  $r_0\in[ \frac{7}{3},\frac{8}{3}]$. We next define  $f$ through    its second derivative.
Set   \begin{align}
   h(r)=& D_0^2e^{D_0(r-r_0)},\ \ \ \ \ r< r_0-D_0^{-1}.\nonumber
   \end{align}
 Let $b\in[ \frac{1}{100},\frac{1}{4}]$ which will be fixed  later. Define further
   \begin{align}
     h(r)=& e^{-1}D_0^2-(1+s)c^{-1}e^{-1}D_0^2(r-r_0+D_0^{-1}),\ \ \ \ \ \ \ \ \ \ \ r\in[r_0-D_0^{-1},r_0-D_0^{-1}+c),\nonumber\\
      h(r)=& - se^{-1}D_0^2 ,\ \ \ \ \ \ \ \ \ \ \ \ \ \ \ \ \ \ \ \
      \ \ \ \ \ \ \ \ \ \ \ \ \ \    \ \ \ \ \ \ \ \ \ \ \ \ \ \ \ \ \  r\in[r_0-D_0^{-1}+c,r_0-c),\nonumber\\
        h(r)=& -se^{-1}D_0^2+c^{-1}e^{-1}(sD_0^2+2eb)(r-r_0+c),\ \ \ \ \ \ \ \   r\in[r_0-c,r_0 ],\nonumber
   \end{align}
where $0<c<(2D_0)^{-1}$  and $s\in (0,1)$ are two  parameters  which will be fixed  in what below.
We have
\begin{align}   \lim_{c\downarrow 0  }\lim_{  s\uparrow 1 }\sup_{\frac{1}{100}\leq b\leq \frac{1}{4}}\int_{r_0-D_0^{-1}}^{r_0} h(r) dr=-e^{-1}D_0  \ \ \ and\ \ \
\lim_{s \downarrow 0  }\inf_{\frac{1}{100}\leq b\leq \frac{1}{4}}\int_{r_0-D_0^{-1}}^{r_0} h(r) dr\geq 0.\label{ccon}
  \end{align}
 Notice also that $-e^{-1}D_0 <-e^{-1}D_0+2r_0b<0$. Therefore, by continuity arguments,   $c$ and $s$ can be fixed by choosing  a solution of the   equation
   \begin{align}
     \int_{r_0-D_0^{-1}}^{r_0} h(r) dr= -e^{-1}D_0+2r_0b.\label{fix}
   \end{align}
Here    $c$ is first fixed as a small  constant in $(0,(2D_0)^{-1})$ so that the rest part of the first term in (\ref{ccon})
is smaller than  $-e^{-1}D_0+2r_0b$. Then  $s$ is further fixed as  a linear  function of $b\in[ \frac{1}{100},\frac{1}{4}]$.

 By definition,  there exists some $r_1\in [r_0-D_0^{-1},r_0-D_0^{-1}+c)$ such that    $h(r_1)=0$.
Then, with (\ref{fix}) and  $2r_0>c$, the definition of $h$ shows that  for  $t\in (r_1,r_0]$,
  \begin{align}\label{sdr}
  \int_{r_0-D_0^{-1}}^{t} h (r)   dr
     \geq&
     \int_{r_0-D_0^{-1}}^{r_0} h (r)   dr -2 bc^{-1}\int_{r_0-c}^{r_0}(r-r_0+c) dr \nonumber\\
     = & -e^{-1}D_0+2r_0b-  bc \nonumber\\
     > &-e^{-1}D_0.
   \end{align}

      Set  \begin{align}
   f(r)=\int_{-\infty}^r dt\int_{-\infty}^t h (u)du,\  \ \ \ \ r\leq r_0.\nonumber
   \end{align}
This definition  shows that   (\ref{fr}) holds.
 With (\ref{fix}), it can be verified  that
 \begin{align}
   f'(r_0-D_0^{-1})=e^{-1}D_0,\  \ \ f''(r_0-D_0^{-1})=e^{-1}D_0^2;\ \ \  \ \  \ \ \lim_{r \uparrow r_0  }f'(r)=2r_0b,\  \ \
\lim_{r \uparrow r_0  }f''(r)=2b.\label{fs}
   \end{align}
  Applying  the first equality above and   (\ref{sdr}), we can check  that    $f$ is   a increasing function  on $(-\infty,r_0]$.

By conditions $D_0>4$, $b\in[ \frac{1}{100},\frac{1}{4}]$, $s\in (0,1)$ and the definition of $h$,
 \begin{align}
  |f''(r)|\leq e^{-1}D_0^2,\ \ \ \ r\in [r_0-D_0^{-1},r_0). \label{fd}
   \end{align}
Consequently,  $f(r)\leq e^{D_0(r-r_0)}$ for $r\leq r_0$.  This inequality and the increasing of $f$   on $(-\infty,r_0]$ show that
    \begin{align}
   e^{-1}\leq f(r)\leq 1, \ \ \ \ r\in [r_0-D_0^{-1},r_0].\label{bod}
   \end{align}
 Notice that by definition
$f(r_0)$ is a linear  function of  $b\in [ \frac{1}{100},\frac{1}{4}]$. Applying    (\ref{bod}) with $r=r_0$ and     continuity arguments,   $b$  can be fixed as
the  solution of the equation   $b= r_0^{-2}f(r_0) $ with
  $b\in [ \frac{1}{100},\frac{1}{4}]$.
 Then  the   estimates  of $b$ in (\ref{sece}) follow from (\ref{bod}) with $r=r_0$ and  $r_0\in[ \frac{7}{3},\frac{8}{3}]$. Set
   \begin{align}
   f(r)=& b r^2,\ \ \ \ r> r_0.\nonumber
   \end{align}
 Since $f(r_0)=b r_0^2 $,  $f$ is continuous at $r_0$. The continuity of $f'$ and $f''$ at $r=r_0$ can be verified by the last two equalities of (\ref{fs}), respectively.  Applying (\ref{fd})
and the left hand side of (\ref{bod}), we can check that  the  second  conclusion of (\ref{sece}) holds. \qed\medskip

\begin{lemma}\label{Pj}
Let $t\geq 0$ and  $\mu$ be a $t$-uniformly log-concave probability measure on   $\mathbb{R}^n$.
 Denote by  $\hat{\mu}$    the projection measure of $\mu$ on a subspace of  $\mathbb{R}^n$.  Then   $\hat{\mu}$  is also  a t-uniformly log-concave probability measure.
\end{lemma}\noindent\textbf{Proof}\  Let $e^{- \varphi}$ be the probability density function of $\mu$.
Let $1\leq k \leq n-1$ and suppose, without loss of generality, that  $\hat{\mu}$ is given by
\begin{align}
d\hat{\mu}(z) =\big(\int_{\mathbb{R}^{n-k}} e^{-\varphi(z,x_{k+1},\cdots
,x_n)}dx_{k+1}\cdots dx_n\big)dz,\ \ \ \
z \in\mathbb{R}^k.\nonumber
 \end{align}

By the  assumption of $\mu$,  there exists some convex function
$V(x)$ such that
\begin{align}
 \varphi(x)= t|x|^2+V(x),\ \ \ \ \ x\in\mathbb{R}^n .\nonumber
 \end{align}
 Then
 \begin{align}
  d\hat{\mu}(z)=\Big(
 \int_{\mathbb{R}^{n-k}}    e^{-t\sum_{i=k+1}^nx_i^2- V(z,x_{k+1},\cdots, x_n )}dx_{k+1}\cdots dx_n \Big)e^{-t|z|^2} dz, \ \ \ \ z \in\mathbb{R}^k. \nonumber
 \end{align}
This equality and  the Pr$\mathrm{\acute{e}}$kopa-Leindler inequality (e.g. Theorem 1.2.3 in \cite{BGVV14}) imply  the conclusion.
\qed\medskip

 \begin{lemma}\label{Step4}
Let   $\mu$ be an isotropic  log-concave probability measure on   $\mathbb{R}^n$.
Let $D_0>4 $,    $r_0\in [ \frac{7}{3},\frac{8}{3}]$ and define a function $f$   by  Lemma \ref{Step5}.
Set
\begin{align}
F_t=\sum_{i=1}^n f(\lambda_{i;t})  ,\ \ \ \ t\geq 0,\label{ft}
\end{align}
where  $0\leq \lambda_{1;t} \leq  \cdots \leq  \lambda_{n;t}$
 are the eigenvalues of the covariance matrix $A_t$, repeated according to their multiplicity.
   Then, for $t_0\in (0,1]$,
   \begin{align}
   \label{in2}
    \mathbb{E} F_t  \leq (t/t_0)^{10^3}   \mathbb{E} F_{t_0} ,\ \ \ \ \ t\in [t_0,t_0\vee D_0^{-4}].
   \end{align}
\end{lemma}\noindent\textbf{Proof}\
 $1)\ \mathbf{a \ known \  differential\ formlula} $

In what below, we need a known differential formula to estimate  $  \mathbb{E} F_t$.
  We   refer to \cite{K21} and \cite{KL22}  for more details and comments  on  it.
 Let    $(e_{i;t})_{1\leq i \leq n}$ be a corresponding orthonormal basis of
   eigenvectors of $A_t$. When the eigenvalues are distinct, taking $\lambda_{i;t}$ as functions of $A_t$,  we have by It$\mathrm{\hat{o}}$ formula
\begin{align}
 d\lambda_{i;t}=\sum_{k=1}^n  u_{iik;t} d B_{k;t}  -\lambda_{i;t}^2dt+\sum_{j:j\neq i}\sum_{k=1}^n\frac{|u_{ijk;t}|^2}{\lambda_{i;t}-\lambda_{j;t}}dt,\label{dif}
 \end{align}
where
\begin{align}
 u_{ijk;t}=\int  \langle x-a_{t},e_{i;t}\rangle \langle x-a_{t},e_{j;t}\rangle\langle x-a_{t},e_{k;t}\rangle p_t(x)dx.\nonumber
 \end{align}
 Notice that $u_{ijk;t}$ is symmetry for indexes  $i,j,k$.
Applying  It$\mathrm{\hat{o}}$ formula with (\ref{dif}), we get
\begin{align}
  dF_t
 =&\sum_{i=1}^n\sum_{k=1}^n u_{iik;t} f'(\lambda_{i;t}) dB_{k,t}-
 \sum_{i=1}^n \lambda_{i;t}^2 f'(\lambda_{i;t})dt+\frac{1}{2}\sum_{i=1}^n\sum_{j\neq i}\sum_{k=1}^n
  |u_{ijk;t}|^2\frac{f'(\lambda_{i;t})-f'(\lambda_{j;t})}{\lambda_{i;t}-\lambda_{j;t}}dt\nonumber\\
 +&\frac{1}{2}\sum_{i=1}^n\sum_{k=1}^n |u_{iik;t}|^2 f''(\lambda_{i;t}) dt .\label{for}
 \end{align}
 The formula above can also  be verified without the distinct assumption of the eigenvalues.
 In the general cases, the quotient term above is interpreted by continuity as $f''(\lambda_{i;t})$ when $\lambda_{i;t}=\lambda_{j;t}$.
 For convenience, this convention will be used  in what below even in the case $i=j$.
Therefore we arrive at the following form
 \begin{align} \label{a1}
  dF_t
  = &\sum_{i=1}^n\sum_{k=1}^n  u_{iik;t}f'(\lambda_{i;t}) dB_{k,t}-
 \sum_{i=1}^n \lambda_{i;t}^2f'(\lambda_{i;t})dt
  +\frac{1}{2}\sum_{i=1}^n\sum_{j=1}^n\sum_{k=1}^n
|u_{ijk;t}|^2\frac{f'(\lambda_{i;t})-f'(\lambda_{j;t})}{\lambda_{i;t}-\lambda_{j;t}}dt\nonumber\\
   :=& \sum_{i=1}^n\sum_{k=1}^n u_{iik;t} f'(\lambda_{i;t}) dB_{k,t}-
 \sum_{i=1}^n \lambda_{i;t}^2f'(\lambda_{i;t})dt +M_tdt.
 \end{align}
 \medskip

 $2)\ \mathbf{ a\ upper\ bound  \ of}$ $ M_t$\ \

 To estimate $M_t$, we make a decomposition first.  By definition, $f'(r)\geq 0$ for $r\geq 0$ and $f'(r)=2br,f''(r)=2b$ for $r\geq r_0$. Then,
\begin{align}\label{dec,01}
&M_t\nonumber\\
= & \frac{1}{2}\sum_{k=1}^n \sum_{ \lambda_{i;t}\leq r_0}  \sum_{\lambda_{j;t}\leq r_0}
|u_{ijk;t}|^2\frac{f'(\lambda_{i;t})-f'(\lambda_{j;t})}{\lambda_{i;t}-\lambda_{j;t}} + \sum_{k=1}^n  \sum_{ \lambda_{i;t}\leq r_0 }  \sum_{  \lambda_{j;t}> r_0}
|u_{ijk;t}|^2\frac{f'(\lambda_{i;t})-f'(\lambda_{j;t})}{\lambda_{i;t}-\lambda_{j;t}}\nonumber\\
+& \frac{1}{2}\sum_{k=1}^n  \sum_{ \lambda_{i;t}> r_0 }  \sum_{  \lambda_{j;t}> r_0}
|u_{ijk;t}|^2\frac{f'(\lambda_{i;t})-f'(\lambda_{j;t})}{\lambda_{i;t}-\lambda_{j;t}}\nonumber\\
\leq &  \frac{1}{2}\sum_{k=1}^n \sum_{ \lambda_{i;t}\leq r_0+ \frac{1}{3} }  \sum_{\lambda_{j;t}\leq r_0+ \frac{1}{3}}
|u_{ijk;t}|^2\Big|\frac{f'(\lambda_{i;t})-f'(\lambda_{j;t})}{\lambda_{i;t}-\lambda_{j;t}}\Big| \nonumber\\
+& \sum_{k=1}^n  \sum_{ \lambda_{i;t}\leq r_0 }  \sum_{  \lambda_{j;t}> r_0+ \frac{1}{3}}
|u_{ijk;t}|^2\frac{f'(\lambda_{i;t})-f'(\lambda_{j;t})}{\lambda_{i;t}-\lambda_{j;t}}
+b\sum_{k=1}^n  \sum_{ \lambda_{i;t}> r_0 }  \sum_{  \lambda_{j;t}> r_0}
   |u_{ijk;t}|^2 \nonumber\\
\leq & \frac{1}{2} \sum_{\lambda_{k;t}\leq r_0+ \frac{1}{3}}  \sum_{ \lambda_{i;t}\leq r_0+ \frac{1}{3}}  \sum_{\lambda_{j;t}\leq r_0+ \frac{1}{3}}
|u_{ijk;t}|^2\Big|\frac{f'(\lambda_{i;t})-f'(\lambda_{j;t})}{\lambda_{i;t}-\lambda_{j;t}}\Big|\nonumber\\
+& \frac{1}{2} \sum_{\lambda_{k;t}> r_0+ \frac{1}{3}}  \sum_{ \lambda_{i;t}\leq r_0+ \frac{1}{3}}  \sum_{\lambda_{j;t}\leq r_0+ \frac{1}{3}}
|u_{ijk;t}|^2\Big|\frac{f'(\lambda_{i;t})-f'(\lambda_{j;t})}{\lambda_{i;t}-\lambda_{j;t}}\Big|
\nonumber\\+&2 b  \sum_{\lambda_{k;t}\leq r_0+ \frac{1}{3}} \sum_{ \lambda_{i;t}\leq r_0 }  \sum_{  \lambda_{j;t}> r_0+ \frac{1}{3}}
 |u_{ijk;t}|^2\frac{\lambda_{j;t}}{\lambda_{j;t}-\lambda_{i;t}}\nonumber\\
+ & 2 b\sum_{ \lambda_{k;t}> r_0+ \frac{1}{3}} \sum_{ \lambda_{i;t}\leq r_0 }  \sum_{  \lambda_{j;t}> r_0+ \frac{1}{3}}
|u_{ijk;t}|^2\frac{\lambda_{j;t}}{\lambda_{j;t}-\lambda_{i;t}}\nonumber\\
+&b\sum_{k=1}^n  \sum_{ \lambda_{i;t}> r_0 }  \sum_{  \lambda_{j;t}> r_0}
   |u_{ijk;t}|^2 \nonumber\\
:=&M_{1,t}+M_{2,t}+M_{3,t}+M_{4,t}+M_{5,t}.
\end{align}
 Here $\sum_{\lambda_{i;t}\leq r_0}$, for example, represents  $\sum_{i:1\leq i\leq n,\lambda_{i;t}\leq r_0}$. This convention
  will be used in what  below.

We have by  the estimate of $f''$ in (\ref{sece}) and that   $f$ is increasing
\begin{align} \label{I1}
M_{1,t}
\leq  &  \frac{1}{2}\sum_{\lambda_{k;t}\leq r_0+ \frac{1}{3}}  \sum_{ \lambda_{i;t}\leq\lambda_{k;t} }  \sum_{\lambda_{j;t}\leq \lambda_{k;t}}
|u_{ijk;t}|^2\Big|\frac{f'(\lambda_{i;t})-f'(\lambda_{j;t})}{\lambda_{i;t}-\lambda_{j;t}}\Big| \nonumber\\
+&  \frac{1}{2}\sum_{\lambda_{i;t}\leq r_0+ \frac{1}{3}}  \sum_{ \lambda_{k;t}\leq\lambda_{i;t} }  \sum_{\lambda_{j;t}\leq \lambda_{i;t}}
|u_{ijk;t}|^2\Big|\frac{f'(\lambda_{i;t})-f'(\lambda_{j;t})}{\lambda_{i;t}-\lambda_{j;t}}\Big| \nonumber\\
+&\frac{1}{2} \sum_{\lambda_{j;t}\leq r_0+ \frac{1}{3}}  \sum_{ \lambda_{k;t}\leq\lambda_{j;t} }  \sum_{\lambda_{i;t}\leq \lambda_{j;t} }
|u_{ijk;t}|^2\Big|\frac{f'(\lambda_{i;t})-f'(\lambda_{j;t})}{\lambda_{i;t}-\lambda_{j;t}}\Big| \nonumber\\
\leq  & \frac{1}{2} D_0^2\sum_{\lambda_{k;t}\leq r_0+ \frac{1}{3}} f(\lambda_{k;t}) \sum_{ \lambda_{i;t}\leq\lambda_{k;t} }  \sum_{\lambda_{j;t}\leq \lambda_{k;t}}
  |u_{ijk;t}|^2  \nonumber\\
+&  \frac{1}{2}D_0^2\sum_{\lambda_{i;t}\leq r_0+ \frac{1}{3}} f(\lambda_{i;t})  \sum_{ \lambda_{k;t}\leq \lambda_{i;t} }  \sum_{\lambda_{j;t}\leq \lambda_{i;t} }
|u_{ijk;t}|^2 \nonumber\\
+& \frac{1}{2}D_0^2 \sum_{\lambda_{j;t}\leq r_0+ \frac{1}{3}} f(\lambda_{j;t})  \sum_{ \lambda_{k;t}\leq \lambda_{j;t} }  \sum_{\lambda_{i;t}\leq \lambda_{j;t}}
|u_{ijk;t}|^2 \nonumber\\
= &\frac{3}{2}D_0^2 \sum_{\lambda_{i;t}\leq r_0+ \frac{1}{3}} f(\lambda_{i;t})  \sum_{ \lambda_{j;t}\leq\lambda_{i;t} }  \sum_{\lambda_{k;t}\leq \lambda_{i;t} }
|u_{ijk;t}|^2.
\end{align}
Similarly, we have by   $r_0+\frac{1}{3}\leq 3$ and $f(3)\leq 2$
\begin{align}\label{I12}
M_{2,t}\leq &\frac{1}{2}D_0^2\sum_{\lambda_{k;t}> r_0+ \frac{1}{3}}  \sum_{ \lambda_{i;t}\leq r_0+ \frac{1}{3} }  \sum_{\lambda_{j;t}\leq r_0+ \frac{1}{3}} f(3)  |u_{ijk;t}|^2\nonumber\\
\leq& D_0^2\sum_{\lambda_{k;t}> r_0+ \frac{1}{3}}  \sum_{ \lambda_{i;t}\leq r_0+ \frac{1}{3} }  \sum_{\lambda_{j;t}\leq r_0+ \frac{1}{3}}  |u_{ijk;t}|^2 .
\end{align}
Since $r_0 \leq \frac{8}{3}$,  we also have
\begin{align}\label{dec31}
M_{3,t}
\leq &18b\sum_{\lambda_{j;t}> r_0+ \frac{1}{3}}\sum_{ \lambda_{i;t}\leq r_0+ \frac{1}{3}} \sum_{ \lambda_{k;t}\leq r_0+ \frac{1}{3}}
 |u_{ijk;t}|^2.
\end{align}
Similarly,
\begin{align}\label{dec1}
M_{4,t}
\leq &\sum_{\lambda_{j;t}> r_0+ \frac{1}{3}}\sum_{ \lambda_{i;t}\leq r_0} \sum_{ \lambda_{k;t}\leq \lambda_{j;t}}
|u_{ijk;t}|^2\frac{\lambda_{j;t}}{\lambda_{j;t}-\lambda_{i;t}}
+\sum_{\lambda_{k;t}> r_0+ \frac{1}{3}} \sum_{\lambda_{i;t}\leq r_0}\sum_{
r_0+ \frac{1}{3}< \lambda_{j;t} \leq \lambda_{k;t}}|u_{ijk;t}|^2\frac{\lambda_{j;t}}{\lambda_{j;t}-\lambda_{i;t}}\nonumber\\
\leq &18 b\sum_{\lambda_{j;t}> r_0+ \frac{1}{3}}\sum_{ \lambda_{i;t}\leq r_0} \sum_{ \lambda_{k;t}\leq \lambda_{j;t}}
 |u_{ijk;t}|^2
+18b\sum_{\lambda_{k;t}> r_0+ \frac{1}{3}} \sum_{\lambda_{i;t}\leq r_0}\sum_{  \lambda_{j;t}\leq \lambda_{k;t}} |u_{ijk;t}|^2
\nonumber\\
\leq &36b\sum_{\lambda_{j;t}> r_0+ \frac{1}{3}}\sum_{ \lambda_{i;t}\leq r_0} \sum_{ \lambda_{k;t}\leq \lambda_{j;t}}
 |u_{ijk;t}|^2 .
\end{align}
We   have
\begin{align}
M_{5,t}\leq&b\sum_{ \lambda_{i;t}> r_0 } \sum_{ \lambda_{k;t}\leq  \lambda_{i;t}}  \sum_{  \lambda_{j;t}\leq  \lambda_{i;t}}
   |u_{ijk;t}|^2 + b\sum_{ \lambda_{j;t}> r_0 } \sum_{ \lambda_{k;t}\leq  \lambda_{j;t}}  \sum_{  \lambda_{i;t}\leq \lambda_{j;t}}
   |u_{ijk;t}|^2\nonumber\\
+ &b\sum_{ \lambda_{k;t}> r_0 } \sum_{ \lambda_{i;t} \leq  \lambda_{k;t}}  \sum_{  \lambda_{j;t}\leq   \lambda_{k;t}}
   |u_{ijk;t}|^2
   \nonumber\\
\leq  &3b\sum_{ \lambda_{i;t}> r_0 } \sum_{ \lambda_{j;t}\leq  \lambda_{i;t}}  \sum_{  \lambda_{k;t}\leq   \lambda_{i;t}}
   |u_{ijk;t}|^2.\label{I5}
\end{align}

$3)\ \mathbf{moments\ estimates}$\ \

 To give further bounds , we prepare a   moment  estimate first. The calculation below relies on  a    technique of moment estimate   used  in, e.g.,  \cite{C21}, \cite{KL22} and  \cite{K23}.
We also need  (\ref{kla}) instead of (\ref{Lic}).   Let $r>0$.
For each $1\leq k \leq n$, set
\begin{align}
g_{r;k}(x)= \sum_{\lambda_{i;t}\leq r}\sum_{  \lambda_{j;t}\leq r}  u_{ijk;t}\langle x-a_{t},e_{i;t}\rangle \langle x-a_{t},e_{j;t} \rangle. \nonumber
\end{align}
By Cauchy-Schwarz inequality,
\begin{align}
   \sum_{  \lambda_{i;t}\leq r} \sum_{ \lambda_{j;t}\leq r}
|u_{ijk;t}|^2 =& \int \langle
x-a_t,e_{k;t}\rangle g_{r;k}(x) p_t(x)dx
\nonumber\\
\leq & \Big(\int\langle
x-a_t,e_{k;t}\rangle^2p_t(x)dx\Big)^{1/2}\big(\mathrm{Var}_{p_t}(g_{r;k})\big)^{1/2}
\nonumber\\
= &\lambda_{k;t}^{1/2}\big(\mathrm{Var}_{ {p}_t}( {g}_{r;k})\big)^{1/2}.\nonumber
 \end{align}
 Let $H$ be the subspace of $\mathbb{R}^n$ spanned by vectors $\{e_{i;t}:{\lambda_{i;t}\leq r}\}$ and assume that it is not an empty set. Denote by   $\hat{g}_{r;k}$ the function of $g_{r;k}$ restricted  on $H$
 and by $\hat{p}_t$ the projection measure of $p_t$ on $H$, respectively. Let $\hat{A}_t$ be the covariance matrix of $\hat{p}_t$. We see that $\|\hat{A}_t\|_{op}\leq r$.
 Applying   (\ref{kla})   with the help of Lemma \ref{Pj}, we further get
\begin{align}
 \sum_{  \lambda_{i;t}\leq r} \sum_{ \lambda_{j;t}\leq r}
|u_{ijk;t}|^2 \leq &\lambda_{k;t}^{1/2}  \big(\mathrm{Var}_{\hat{p}_t}(\hat{g}_{r;k})\big)^{1/2}\nonumber\\
\leq  &\lambda_{k;t}^{1/2}\Big(t^{-1/2}r^{1/2} \int_{H} |\nabla \hat{g}_{r;k}(\hat{x})|^2\hat{p}_t(\hat{x})d\hat{x} \Big)^{1/2}\nonumber\\
=  &t^{-1/4}r^{1/4}\lambda_{k;t}^{1/2}\Big( \int  |\nabla  {g}_{r;k}(x)|^2 {p}_t(x)dx \Big)^{1/2}\nonumber\\
\leq  &t^{-1/4}r^{1/4}\lambda_{k;t}^{1/2} \Big(4 r \sum_{
\lambda_{i;t}\leq r} \sum_{ \lambda_{j;t}\leq r } |u_{ijk;t}|^2
\Big)^{1/2},\nonumber
\end{align}
which implies that
\begin{align}\label{ba1}
 & \sum_{  \lambda_{i;t}\leq r} \sum_{ \lambda_{j;t}\leq r}
|u_{ijk;t}|^2  \leq  4t^{-1/2} r^{3/2}\lambda_{k;t}.
\end{align}

$4)\ \mathbf{combing  \ the\ three\ parts\ above}$\ \

It follows from  (\ref{I1}), (\ref{ba1})  and $r_0\leq \frac{8}{3}$ that \begin{align}
M_{1,t}
    \leq&  54 \sqrt{3}D_0^2t^{-1/2}  \sum_{\lambda_{i;t}\leq r_0+ \frac{1}{3}}f(\lambda_{i;t}).\nonumber
\end{align}
Applying (\ref{I12}), (\ref{ba1}),  the lower bound of $b$ in  (\ref{sece}) and $r_0\geq \frac{7}{3}$, we have \begin{align}
M_{2,t}
  \leq& 12 \sqrt{3}D_0^2 t^{-1/2}\sum_{\lambda_{k;t}> r_0+ \frac{1}{3}} \lambda_{k;t}
  \leq 90 \sqrt{3} D_0^2t^{-1/2}  \sum_{\lambda_{i;t}> r_0+ \frac{1}{3}}f(\lambda_{i;t}).\nonumber
\end{align}
Similarly,   by  (\ref{dec31}) and (\ref{ba1}),
\begin{align}
M_{3,t}\leq &  216\sqrt{3}t^{-1/2}\sum_{\lambda_{j;t}> r_0+ \frac{1}{3}} b\lambda_{j;t}
  \leq 81 \sqrt{3}t^{-1/2} \sum_{\lambda_{i;t}> r_0+ \frac{1}{3}}f(\lambda_{i;t}).\nonumber
\end{align}
Applying (\ref{dec1}), (\ref{ba1}) and (\ref{Lic}), we get
\begin{align}
&M_{4,t}\leq  144 t^{-1/2}   \sum_{\lambda_{j;t}> r_0+ \frac{1}{3} }b \lambda_{j;t}^{5/2}
\leq 144t^{-1}\sum_{\lambda_{i;t}> r_0+ \frac{1}{3}}f(\lambda_{i;t}),\nonumber
\end{align}
and get with (\ref{I5})
\begin{align}
&M_{5,t}
\leq 12 t^{-1}\sum_{\lambda_{i;t}> r_0}f(\lambda_{i;t}).\nonumber
\end{align}

Let $t_0\in (0,1]$. Combing   (\ref{a1}), (\ref{dec,01}) and the estimates of
$M_{1,t},\cdots,M_{5,t}$ stated  above, we obtain that for $t\in [t_0, 1]$,
 \begin{align}
 \frac{d}{dt}  \mathbb{E} F_t\leq & (156+81\sqrt{3t})t^{-1} \sum_{i=1}^n \mathbb{ E}f(\lambda_{i;t})+ 90\sqrt{3}D_0^2 t^{-1/2} \sum_{i=1}^n \mathbb{ E}f(\lambda_{i;t}),\label{tb}\\
 \leq &\big(500+90\sqrt{3}D_0^2 t^{-1/2}\big) \mathbb{E} F_t\nonumber\end{align}
which further gives for   $t\in [t_0, t_0\vee D_0^{-4}]$,
 \begin{align}
 \frac{d}{dt}  \mathbb{E} F_t\leq &
 10^3 t^{-1} \mathbb{E} F_t.\nonumber
\end{align}
Then we get (\ref{in2}) by the estimate above.
 \qed

\section{Proof of Lemma \ref{Step0}   }

$1)\ \mathbf{a\ known\ bound\ for \  the \ operator\ norm}$\ \

The following  parameter  is from \cite{E13}:
\begin{align}
\kappa_n^2=\sup_X   \|\mathbb{E}X_1(X\otimes X)\|_2^2, \label{kappa}\end{align}
where $X$ runs over all  isotropic log-concave random variables $X$
in $\mathbb{R}^n$. By the  inequality    (3.14) in  \cite{K23} and Theorem 1.2   in  \cite{K23},
there exists some universal constant $C_1>0$ such that
\begin{align}
 \kappa_n^2\leq C_1 \log n,\ \ \ n\geq 2.\nonumber
\end{align}By Lemma 5.2 in \cite{KL22} and the bound of $\kappa_n$ above, there
exists some universal constant $C_2>0$ such that
\begin{align}
 \mathbb{P}(\|A_t\|_{op}\geq 2)\leq e^{-(C_2t)^{-1}}, \ \ \ \ \ t\in[0,t_1],\label{oper}
\end{align}
where $t_1= C_2^{-1}(\log n)^{-2} $ with $n\geq 2$.

In what below we assume that $n$ is bigger than a universal constant $C$ such that
\begin{align}
   t_1<&e^{-10^3},  \label{t1}\\
  n\geq&    e^{10|\log t_1|^2} , \label{nc}
\end{align}
Therefore,  $t_1^{-2}e^{-(C_2t_1)^{-1}}\leq 1$. By differentiation, we further get
\begin{align}
 t^{-2}e^{-(C_2t)^{-1}}\leq 1,\ \ \ \ \ \ t\in(0,t_1]. \label{ncn}
\end{align}

$2)\ \mathbf{  assistant\ stochastic\ processes }$\ \

With $  t_1=C_2^{-1}(\log n)^{-2}$ as above, set
\begin{align}
\  \ \ \   t_{k+1}=|\log t_{k}|^{-16},\ \ \ \ \ \ k\geq 1.\label{ttkk}
   \end{align}
 Define
      \begin{align}
  k_0=\sup \{k\geq 1: t_k\leq e^{-10^3}\},\nonumber
   \end{align}
and   set
   \begin{align}
  s_1=\frac{7}{3}, \ \ \   s_{k+1}= \frac{7}{3}+ \sum_{i=2}^{k+1}|\log t_{i}|^{-1/2}  ,\ \ \ \ \ \ k= 1,\cdots,k_0.\nonumber
   \end{align}

Noticing  that  $h(r)= r-32\ln r$ is a increasing function of $r$ in  $(32,+\infty)$ and  $t_{k_0+1}\leq e^{-100}$,
 we have   $-t_{k+1}^{-\frac{1}{16}}\leq -32\ln t_{k+1}^{-\frac{1}{16}}$ when $1\leq k\leq k_0$.   By the relation $t_k=e^{-t_{k+1}^{-\frac{1}{16}}}$, we further get
   \begin{align}
 t_k\leq (t_{k+1})^2  ,\ \ \ \ \ \ k= 1,\cdots,k_0.\label{tk}
   \end{align}
  Consequently, for $1\leq k\leq k_0$ with $k_0\geq 2$,
      \begin{align}
  \frac{7}{3}\leq s_k\leq  s_{k_0}=\frac{7}{3}+ \sum_{k=2}^{k_0}(t_{k+1})^{\frac{1}{32}}\leq
   \frac{7}{3}+ \sum_{k=2}^{k_0}(t_{k_0+1})^{ 2^{k_0-k-5}} \leq \frac{7}{3}+ \sum_{k=2}^{k_0} e^{-100\cdot 2^{k_0-k-5}}\leq
 \frac{8}{3}.\label{inc}
   \end{align}
The estimates above hold also when  $k_0=1$.

In what below we continue to use the notations in Lemma \ref{Pj}. With (\ref{inc}) we can set
   \begin{align}
  F_{k;t}= \sum_{i=1}^nf_{k}(\lambda_{i;t}) ,\ \ \ \ \ t\geq 0,\ 1\leq k\leq k_0,\label{ap}
   \end{align}
   where $f_k$ is a function     in Lemma \ref{Step5} with $r_0=s_k$ and $D_0=|\log t_k|^4$.
Notice that $D_0^{-4}=t_{k+1}$ by (\ref{ttkk}).\medskip

$3)\ \mathbf{estimates\ for\ the\ expectation\  of} $   $(F_{k;t})$ $  \mathbf{ by \ finite \ induction}$\ \

It follows from  (\ref{tk}) and $t_{k_0+1}\leq e^{-100}$ that  $t_k<t_{k+1}$ for
$k= 1,\cdots,k_0$. By definition, $ e^{(C_2t_1)^{-1}} \geq n$. Then, for $t\in [t_1,t_2]$, we have by Lemma
\ref{Step4} with $f=f_1$, $t_0=t_1$, Lemma \ref{Step5} and  (\ref{oper})
  \begin{align}
    \mathbb{E} F_{1,t}  \leq   & t_1^{-10^3}  \mathbb{E} F_{1,t_1} \nonumber\\
    = & t_1^{-10^3} \mathbb{E} F_{1,t_1} I_{\|A_{t_1}\|_{op}\leq 2}+   t_1^{-10^3} \mathbb{E} F_{1,t_1} I_{\|A_{t_1}\|_{op}> 2}\nonumber\\
\leq &t_1^{-10^3} \sum_{i=1}^n  \mathbb{E}  f_1(\lambda_{i;t_1}) I_{\lambda_{i;t_1}\leq 2 }+t_1^{-10^3-2} n   \mathbb{P}(\|A_{t_1}\|_{op}\geq 2)\nonumber\\
\leq &t_1^{-10^3}     f_1(2)n+t_1^{-10^3-2}  \nonumber\\
    = &  t_1^{-10^3}e^{-\frac{1}{3}|\log t_1|^4}n+ t_1^{-10^3-2}  \nonumber,
   \end{align}
   where   (\ref{Lic}) is used  in the second  inequality above.
Applying    (\ref{t1}) and (\ref{nc}), we further get
  \begin{align}
    \mathbb{E} F_{1,t}  \leq
      &  e^{-\frac{1}{3}|\log t_1|^4+10^3|\log t_1|}n+ e^{-10|\log t_1|^2+(10^3+2)|\log t_1|}n\nonumber\\
    \leq &  e^{-|\log t_1|^2}n.\label{0}
   \end{align}
 We claim further that      for every integer $k\in \{1,\cdots,k_0\}$,
     \begin{align}
       \mathbb{E} F_{k,t}
    \leq &  e^{-|\log t_{k}|^2}n, \ \ \ \ \ \ t\in [t_{k},t_{k+1}].\label{intr}
   \end{align}  By (\ref{0}),
   the estimates above hold when $k=1$.
Suppose next that (\ref{intr})    holds for some  $k\in \{1,\cdots,k_0-1\}$. To verify  the assertion  above we only
need to show  that (\ref{intr}) still holds when  $k$ is replaced by
$k+1$.

Recall that $f_k$ is defined by Lemma \ref{Step5}  with $r_0=s_k$, $s_k\in [\frac{7}{3},\frac{8}{3}]$. Denote the corresponding  parameter  $b$ of $f_k$ by $b_k$ in what below.
Since $f_{k+1}$ is a increasing function,  for $r\leq s_{k+1}$, we have $f_{k+1}(r)\leq f_{k+1}(s_{k+1})=b_{k+1}s_{k+1}^2$.
Then it follows   from the definitions of $f_{k+1}$ and $f_k$,    the   estimates
of $b$ in  (\ref{sece}) and (\ref{inc}) that
 \begin{align}
    & \sum_{i=1}^n  \mathbb{E}  f_{k+1}(\lambda_{i;t_{k+1}}) I_{\lambda_{i;t_{k+1}}> s_{k}}\nonumber\\
    =& \sum_{i=1}^n  \mathbb{E}  f_{k+1}(\lambda_{i;t_{k+1}}) I_{ s_{k}<\lambda_{i;t_{k+1}}\leq s_{k+1}}+ \sum_{i=1}^n  \mathbb{E}  f_{k+1}(\lambda_{i;t_{k+1}}) I_{\lambda_{i;t_{k+1}}> s_{k+1}}\nonumber\\
             \leq&    \sum_{i=1}^n \mathbb{E} b_{k+1}s_{k+1}^2I_{s_k<\lambda_{i;t_{k+1}}\leq  s_{k+1}}+  \sum_{i=1}^n \mathbb{E}  b_{k+1}\lambda_{i;t_{k+1}}^2I_{\lambda_{i;t_{k+1}}> s_{k+1}}\nonumber\\
          \leq&  6\sum_{i=1}^n \mathbb{E}  b_{k}\lambda_{i;t_{k+1}}^2I_{s_k<\lambda_{i;t_{k+1}}\leq  s_{k+1}}+4 \sum_{i=1}^n \mathbb{E}  b_{k}\lambda_{i;t_{k+1}}^2I_{\lambda_{i;t_{k+1}}> s_{k+1}}\nonumber\\
          \leq &6 \mathbb{E} F_{k,t_{k+1} } ,\nonumber
   \end{align}
where the increasing of $f_{k+1}$ is used in the first inequality. Here $I_{\lambda_{i;t_{k+1}}> s_{k}}$  refers to  the indicator function of the set  $\{i:\lambda_{i;t_{k+1}}> s_{k},1\leq i\leq n\}$ and this convention is adopted in the whole proof.

Let $t\in [t_{k+1},t_{k+2}]$.   It follows by    Lemma \ref{Step4} with $f=f_{k+1}$,
$t_0=t_{k+1}$, the definition of $f_{k+1}$ and the estimate above that
 \begin{align}
   & \mathbb{E} F_{k+1,t}\nonumber\\\leq &
     (t_{k+1})^{-10^3} \mathbb{E} F_{k+1,t_{k+1}} \nonumber\\
      =&(t_{k+1})^{-10^3} \sum_{i=1}^n  \mathbb{E} f_{k+1}(\lambda_{i;t_{k+1}}) I_{\lambda_{i;t_{k+1}}\leq s_{k}}+ (t_{k+1})^{-10^3} \sum_{i=1}^n
       \mathbb{E}  f_{k+1}(\lambda_{i;t_{k+1}}) I_{\lambda_{i;t_{k+1}}> s_{k}}\nonumber\\
     \leq &(t_{k+1})^{-10^3}  f_{k+1}(  s_{k})n+ +6(t_{k+1})^{-10^3} \mathbb{E} F_{k,t_{k+1} }\nonumber\\
         =&   (t_{k+1})^{-10^3}  e^{-|\log t_{k+1}|^{\frac{7}{2}}}n +6(t_{k+1})^{-10^3} \mathbb{E} F_{k,t_{k+1} }
\nonumber\\
         \leq&    \frac{1}{2} e^{-|\log t_{k+1}|^2}n +6(t_{k+1})^{-10^3} \mathbb{E} F_{k,t_{k+1} },\label{k+1}
   \end{align}
where $t_{k+1}\leq e^{-10^3}$ is used in the last step above.

It can be verified that  $h(r)= \frac{1}{2}r^{-1}+8\cdot 10^3\ln r$ is a decreasing function on $(0,e^{-10^2})$. Since $t_{k_0}\leq e^{-10^3}$,  this yields
 $\frac{1}{2}(t_{k+1})^{-\frac{1}{8}}-8\cdot 10^3\ln t_{k+1}^{-\frac{1}{8}} \geq 0$.
Noticing  that the function $h(r)= \frac{1}{\sqrt{2}}r^{-\frac{1}{16}}+\ln r$ is decreasing on $(0,e^{-10^3})$, we can also check that
 $\frac{1}{2}(t_{k+1})^{-\frac{1}{8}}-|\ln t_{k+1}|^2 \geq 3$.  Therefore,
the induction assumption of (\ref{intr})    and (\ref{k+1})  give \begin{align}
    \mathbb{E} F_{k+1,t}
     \leq&   \frac{1}{2} e^{-|\log t_{k+1}|^2}n+6(t_{k+1})^{-10^3} e^{-|\log t_k|^2}n \nonumber\\
     = & \frac{1}{2} e^{-|\log t_{k+1}|^2}n+6  e^{-(t_{k+1})^{-\frac{1}{8}}+10^3|\log t_{k+1}|}n\nonumber\\
        \leq & \frac{1}{2}  e^{-|\log t_{k+1}|^2}n+  6e^{-\frac{1}{2}(t_{k+1})^{-\frac{1}{8}}  }n  \nonumber\\
     \leq & e^{-|\log t_{k+1}|^2}n,\nonumber
   \end{align}
   which completes the induction.\medskip

$4)\ \mathbf{estimates\ of\  }\mathbb{E}\mathrm{Tr}( A_t^2)$\ \

By definition, for $t\in [t_1,e^{-10^{3}}]$, there exists some $  k \in\{1,\cdots, k_0\}$ such that
  $t\in [t_k,t_{k+1})$. Then  the lower bound of $b$ in    (\ref{sece}),  (\ref{intr})  and $t_{k_0}\leq e^{-10^3}$yield
\begin{align}
 \mathbb{E}\mathrm{Tr}( A_t^2)
 = & \sum_{ \lambda_{i;t}\leq \frac{8}{3}}  \mathbb{E}  \lambda_{i;t}^2+
 \sum_{ \lambda_{i;t}> \frac{8}{3}}  \mathbb{E}  \lambda_{i;t}^2
 \nonumber\\
 \leq &   \frac{64}{9}n+20\mathbb{E}  F_{k,t} \nonumber\\
 \leq &\frac{64}{9}n+20e^{-|\log t_k|^2}n\nonumber\\
 \leq &8n.\label{b1}
\end{align}
Applying  (\ref{Lic}), (\ref{oper}) and  (\ref{ncn}), the bound above holds also for $t\in (0,t_1)$.

To get the lower bound, we use the following formula in  Lemma 29   \cite{LV24}:
\begin{align}
dA_t=   \int(x-a_t)(x-a_t)^T  dp_t(x)dx-A_t^2dt.\label{b2}
\end{align}
Since $A_0$ is equal to the identity,  it follows from  (\ref{b1}) and (\ref{b2}) that
\begin{align}
\mathbb{E}\mathrm{Tr} (A_t)\geq \frac{n}{\sqrt{2}} ,\ \ \ \ t\in [0,e^{-10^3}).\label{idd}
\end{align}
Notice that the  formula (\ref{b2}) is also a ingredient  to deduce (\ref{a1}). The use of (\ref{b2}) in (\ref{idd}) can  be replaced by (\ref{a1}) with $f$ being  the identity function.
Consequently, for $t\in [0,e^{-10^3})$,
\begin{align}
\mathbb{E}\mathrm{Tr}( A_t^2)\geq \sum_{i=1}^n\big(\mathbb{E}  A_{ii;t}\big)^2
\geq n^{-1}\big(\sum_{i=1}^n\mathbb{E}  A_{ii;t}\big)^2=
 n^{-1}\big(\mathbb{E}\mathrm{Tr} A_t\big)^2\geq \frac{n}{2} ,\nonumber
\end{align}
which finish the proof
\qed  \medskip

There is an obstruction to bound   $ \mathbb{E}\|A_t\|_{op}$ by a universal constant when $t>(\log n)^{-1}$. See
\cite{KL24} for more details. The following result is an
upper bound   of this obstruction.
\begin{corollary} \label{comb}
Let $N>1$ and $C_2$ be the constant in (\ref{oper}). Set $\tilde{c}=1\wedge C_2^{-1}$. Then there exists a   constant $C>3$ depending on $N$ such that, under  the same assumptions as Lemma \ref{Step4},  \begin{align}
 \mathbb{E}  \sum_{ \lambda_{i;t}\geq  \frac{8}{3}}  \mathbb{E}  \lambda_{i;t}^2
 \leq & 20 \cdot  e^{-10^{-4}N^{-1/2}(\log n)^{1/2}}n,\ \ \ \ \label{b11}
\end{align}
where $t\in [\tilde{c}(\log n)^{-2},1\wedge (N(\log n)^{-1})]$ and $n\geq C$.

 \end{corollary}\noindent\textbf{Proof}\ Set  $t_0=\tilde{c}(\log n)^{-2}$.
Let $f$ be the function    in Lemma \ref{Step5} with $r_0=\frac{8}{3}$ and $D_0= 2^{-1}10^{-3}N^{-1/2}(\log n)^{1/2}$.
Define $F_t$ by (\ref{ft}). Assume that $n$ is big enough such that $t_0<1$ and  $D_0\geq \frac{3}{2}$.
It follows from the definition of $f$, (\ref{Lic}) and  (\ref{oper})  that \begin{align}
    \mathbb{E} F_{t_0}
      =&  \sum_{i=1}^n  \mathbb{E} f (\lambda_{i;t_0}) I_{\lambda_{i;t_0}\leq 2}+   \sum_{i=1}^n
       \mathbb{E}  f (\lambda_{i;t_0}) I_{\lambda_{i;t_0}>  2}\nonumber\\
\leq &  f (2) n +   P(\|A_{t_0}\|_{op}>  2)\sum_{i=1}^n
\lambda_{i;t_0}^2\nonumber\\\leq
&   e^{-\frac{1}{3}10^{-3}N^{-1/2}(\log n)^{1/2}}n+  \tilde{c}^{-2}(\log n)^4 e^{-\tilde{c}(\log n)^2} n.\label{t11}
   \end{align}

For $t\in [t_0,  1]$, we have by (\ref{tb})
 \begin{align}
 \frac{d}{dt}  \mathbb{E} F_t\leq & (156+81\sqrt{3t})t^{-1}\mathbb{E} F_t+ 90\sqrt{3}\cdot 2^{-2}10^{-6} N^{-1}(\log n)t^{-1/2} \mathbb{E} F_t\nonumber\\
\leq & 10^3t^{-1}\mathbb{E} F_t+ 10^{-4}N^{-1} (\log n) t^{-1/2} \mathbb{E} F_t\nonumber .
\end{align}
  For $t\in [t_0,1\wedge (N(\log n)^{-1})]$, this  implies  that
  \begin{align}
   \mathbb{E}  \sum_{ \lambda_{i;t}\geq  \frac{8}{3}}  \mathbb{E}  \lambda_{i;t}^2
 \leq & 20   \mathbb{E} F_t \nonumber\\
 \leq   &20 (\frac{t}{t_0})^{10^3}   e^{2\cdot 10^{-4}N^{-1} (\log n) t^{1/2}} \mathbb{E} F_{t_0}\nonumber\\
\leq &   20(\tilde{c}^{-1}N\log n)^{10^3} e^{\frac{1}{5}\cdot 10^{-3}N^{-1/2}(\log n)^{1/2}} \mathbb{E} F_{t_0} .\nonumber
   \end{align}
Applying (\ref{t11}) and the estimate above, we get the conclusion by choosing $n$ large enough depending on $N$.
\medskip\qed

\section{$ \log(\log n ) $\ bound \ for \ the \ thin\ cell \ constant\
   }
 }

Let  $s\geq 0$.
 We know   from  \cite{KP21}
 that $C_P(\mu)\leq  C_P(\mu_s)$. For  an isotropic  log-concave probability measure $\mu$ on   $\mathbb{R}^n$, direct calculation shows that
the covariance matrix of $\mu_s$ is $(1+s)Id$. The relation between the thin cell parameter of  $\mu_s$
and $\mu$ is given in (\ref{cn}) below.
 \begin{lemma}\label{Step1A}
Let   $\mu$ be an isotropic  log-concave probability measure on   $\mathbb{R}^n$ with density   $\rho$.
Then for $s\geq 0$,
\begin{align}
  \int  (|x|^2-(1+s)n)^2\rho_s(x)dx
    =&  \int   (|y|^2-n)^2 \rho(y)dy+ 2(s^2+2s)n.\label{cn}
 \end{align}Moreover,
\begin{align}
  -s^{-1}Id\leq\nabla^2\log\rho_s(x)\leq 0,\ \ \ \ \ x\in \mathbb{R}^n,\ \ s>0.\label{Id}
 \end{align}
  Assume further  that     $\rho=e^{-\varphi}$ for some smooth function $\varphi$. Then for any unit vector $\theta\in \mathbb{R}^n$,
\begin{align}
 \int  \partial_\theta^2 \varphi d\mu\geq 1.\label{theta}
 \end{align}

\end{lemma}\noindent\textbf{Proof}\
By definition,
 \begin{align}
  &\int  (|x|^2-(1+s)n)^2\rho_s(x)dx \nonumber\\=
  & \frac{ 1}{(2\pi s)^{n/2}}\int   \rho(y)dy\int (|x|^2-(1+s)n)^2
 \exp\{-\frac{|y-x|^2}{2s}  \}dx\nonumber
  \\=
  &   \frac{ 1}{(2\pi s)^{n/2}}\int   \rho(y)dy\int  \Big((|x- y|^2-sn)+ (|y|^2- n)+2 \langle x- y,y\rangle\Big)^2
 \exp\{-\frac{|x- y|^2}{2s}  \}dx\nonumber
\\
 =& \int  \Big( 2s^2n+(|y|^2-n)^2  +4 s|y|^2\Big)\rho(y)dy
\nonumber\\
 =&  \int   (|y|^2-n)^2 \rho(y)dy+ 2(s^2+2s)n,\nonumber
 \end{align}
which gives the first  conclusion.

Let $\theta$ be a unit vector in $\mathbb{R}^n$. For $x\in \mathbb{R}^n$,
\begin{align}
  \partial_\theta^2\log\rho_s(x)=&
 \frac{ -1}{(2\pi s)^{n/2}\rho_s(x) }\int \big(  s^{-1}- s^{-2}( \langle y,\theta\rangle-m_\theta(x))^2\big) \rho(y)
 \exp\{-\frac{|y-x|^2}{2s}  \}dy\label{par} ,
 \end{align}
where
\begin{align}
m_\theta(x)=&
 \frac{ 1}{(2\pi s)^{n/2}\rho_s(x) }\int  \langle y,\theta\rangle\rho(y)
 \exp\{-\frac{|y-x|^2}{2s}  \}dy\nonumber .
 \end{align}
This expression shows that    the left hand side of (\ref{Id}) holds.
The right hand side of (\ref{Id}) is well known. It is a consequence of (\ref{Lic}) and (\ref{par}) or the Pr$\mathrm{\acute{e}}$kopa-Leindler inequality.

Integration by parts and Cauchy-Schwarz inequality show that
\begin{align}
 \int  \partial_\theta^2 \varphi(x)d\mu=   \int  \big(\partial_\theta  \varphi(x)\big)^2d\mu
\geq  \Big(\int    \langle x, \theta \rangle^2  d\mu \Big)^{-1}\Big(\int   \langle x, \theta \rangle\partial_\theta  \varphi(x)d\mu \Big)^2
=1,\nonumber
 \end{align}
which gives the last conclusion.
\qed\medskip

 \begin{lemma}\label{Tecn}
Let  $h\geq 0$ be a differential and non-increasing function on $[0,\infty)$. Then for any integer $N\geq 0$,
   \begin{align}
   \int_0^{2^N} (s^{-1/2}\wedge 1)\big( - h'(s)\big)^{1/2} ds \leq  h(0)^{1/2} (N+1)^{1/2}.\nonumber
   \end{align}

\end{lemma}\noindent\textbf{Proof}\
It follows from the assumption and  Cauchy-Schwarz inequality that
  \begin{align}
   &\int_0^{2^N} (s^{-1/2}\wedge 1)\big( - h'(s)\big)^{1/2} ds \nonumber\\
   = & \int_0^1\big( - h'(s)\big)^{1/2} ds+ \sum_{k=1}^N \int_{2^{k-1}}^{2^k} s^{-1/2}\big( - h'(s)\big)^{1/2} ds
   \nonumber\\
   \leq &  (h(0)-h(1))^{1/2} +\sum_{k=1}^N 2^{-(k-1)/2} \int_{2^{k-1}}^{2^k}  \big( - h'(s)\big)^{1/2} ds \nonumber\\
     \leq &   (h(0)-h(1))^{1/2}+\sum_{k=1}^N \Big( \int_{2^{k-1}}^{2^k}   - h'(s)  ds\Big)^{1/2} \nonumber\\
   = &  (h(0)-h(1))^{1/2}+ \sum_{k=1}^N  (h(2^{k-1})-h(2^k))^{1/2} \nonumber\\
   \leq &   h(0)^{1/2} (N+1)^{1/2} \nonumber,
   \end{align}
which completes the proof.
\qed \medskip

Recall that $(Q_s)$ is defined by (\ref{st}). To use the $\Gamma$-calculus of $(Q_s)$ developed in  \cite{KP21}, we first restate some definitions and results    in \cite{K23} and \cite{KP21}.
 A log-concave probability measure $\mu$ is called
regular if its density $\rho$ is smooth and positive in $\mathbb{R}^n$ and the following two conditions hold:

i) for some $\varepsilon>0$, $\varepsilon I_d\leq \nabla^2 \varphi\leq \varepsilon^{-1}Id$, where $\varphi=-\log \rho$;

ii) the function $\varphi$ and each of   its partial derivatives of all orders grows at most polynomially at infinity.

We say that a measurable function $u$ has  subexponential decay  relative to $\rho$ if for any $a>0$ there exists some $C>0$ such that
  \begin{align}
  |u(x)|\leq \frac{C}{\sqrt{\rho(x)}} e^{-a|x|},\ \ \ \ \ when\   \rho(x)> 0.\label{se}
\end{align}
A smooth function is called $\mu$-tempered if all of its partial derivatives of all orders    have  subexponential decay  relative to $\rho$.
 Denote
  \begin{align}
  &  \mathcal{F}_\mu=\{u: u\ is \ a \  \mu \textit{{-}}tempered \ function\} .\label{mu}
 \end{align}

 When $\mu$ is regular,  the operator $L$ in (\ref{gen}) can be extended from $\mathcal{F}_\mu$   as    a
self-adjoint operator on $L^2(\mu)$. Moreover, it  has a discrete spectrum and  $\mu$-tempered eigenfunctions. See Lemma 2.2  in  \cite{K23} and the references therein for more details.

Let $H^1(\mu)$ be the space of all the functions in $L^2(\mu)$ whose  weak derivatives belongs to $L^2(\mu)$.
Define the    Rayleigh quotient of $u\in H^1(\mu)$ at time $s\geq 0$  by
    \begin{align}
  &  R_s u =\frac{ \int|\nabla Q_s u |^2 d\mu_s}{ \int|  Q_s u |^2 d\mu_s} .\nonumber
 \end{align}By  a derivative  estimate of  Rayleigh quotient in \cite{KP21}, it shows in Lemma 3.1 \cite{KL22} that,  for   $ u\in H^1(\mu)$ with $ \|  u\|_{L^2(\mu)}=1$,
  \begin{align} & \|Q_s u\|_{L^2(\mu_s)}^2\geq \exp\{-s\int |\nabla u|^2 d\mu\},\ \ \ \ \ s\geq 0. \label{control}
 \end{align}

Let $s\geq 0$. Recall that $\Gamma_{0,s}(u,v)=uv$. For smooth functions $u$ and $v$, it shows in \cite{KP21} that $\Gamma_{1,s}(u,v)=  \langle \nabla u,\nabla v\rangle  $ and
 \begin{align}
   \Gamma_{2,s}( u,u) =\|\nabla^2 u\|_{HS}^2-2 \langle \nabla^2 \log \rho_s \cdot \nabla u, \nabla u\rangle
 ,
   \label{gat}
 \end{align} where $\|\cdot\|_{HS}$ represents the Hilbert-Schmidt norm.
 Define the Laplacian operator $L_s$ of $\mu_s$   by   \begin{align}
L_s =\Delta +\nabla \log \rho_s \cdot\nabla,\ \ \ s> 0.\nonumber
 \end{align}
 The   Bochner formula of $L_s$  in (33) of \cite{KP21} reads
   \begin{align}
 \int (L_s u)^2  d\mu_s=\int\Big(
   \|\nabla^2 u\|_{HS}^2- \langle \nabla^2 \log \rho_s  \cdot \nabla u, \nabla u\rangle\Big)d\mu_s,\ \ \ u\in \mathcal{F}_{\mu_s},\ s>0.
\label{boc}
 \end{align}
For $\lambda>0$, denote by $E_\lambda $ the spectral projection below level $\lambda$ for   $-L$.
 \begin{lemma}\label{TecN} Let   $\mu$ be a regular   log-concave probability measure on   $\mathbb{R}^n$.
 Let $g_k(x)=x_k$ for $x\in \mathbb{R}^n$ and $1\leq k \leq n$.
Assume also that  for some $c_1\in(0,e^{-1})$,
        \begin{align}
   \sum_{k=1}^n\int   |Q_s g_k|^2 d \mu_s \leq c_1^{-1}(s^{-1} \wedge 1)n,\ \ \ \ \ \ \ s> 0 .\label{HEEE}
 \end{align}
 Then   for  $\lambda\in (0, e^{-1})$,
    \begin{align}
   \sum_{k=1}^n\int   |  g_{k,\lambda}|^2 d \mu  \leq  10^2 c_1^{-2}\lambda|\log\lambda|
   n,\ \ \ \ \ \ \ \label{gk2}
 \end{align}
  where
    \begin{align}
   g_{\lambda,k} =E_\lambda g_k,\ \ \ \ \ \ \ 1\leq k\leq n.\nonumber
 \end{align}
\end{lemma}\noindent\textbf{Proof}\
Let  $\lambda\in (0,  e^{-1})$ and assume without loss of generality that the left hand side of (\ref{gk2}) is positive. Define  $c_0=c_0(\lambda) $    by
    \begin{align}
  \sum_{k=1}^n \int| g_{k,\lambda}|^2 d \mu=c_0^{-1}\lambda n.\label{gk}
 \end{align}
Assume further that $c_0<e^{-4}$. Otherwise,  the Lemma holds automatically from $\lambda,c_1\in (0, e^{-1})$.
Since $\mu$ is regular, we know from Lemma 2.2  in  \cite{K23} that  $g_{k,\lambda},g_k-g_{k,\lambda}\in \mathcal{F}_\mu$ for $1\leq k \leq n$.
By Lemma 2.2 in \cite{KP21}, we also have $Q_sg_{k,\lambda} \in \mathcal{F}_{\mu_s}$ for $1\leq k \leq n$ and $s>0$.

Set $s_1= c_0^{1/2} \lambda^{-1}$. Since $E_\lambda$ is a spectral projection, we have by     (\ref{di})   with $k=0$
  \begin{align}
&\Big| \sum_{k=1}^n\int  \big( Q_{s_1}(g_k-g_{k,\lambda})\big)  Q_{s_1}g_{k,\lambda}  d \mu_{s_1}\Big| \nonumber\\
=&\Big| \sum_{k=1}^n\int_0^{s_1}ds \int \langle\nabla Q_s(g_k-g_{k,\lambda}), \nabla Q_sg_{k,\lambda}\rangle d \mu_s\Big|\nonumber\\
\leq&  \Big| \sum_{k=1}^n\int_0^{s_1}ds \int \langle\nabla Q_sg_k , \nabla Q_sg_{k,\lambda}\rangle d \mu_s\Big|
+  \sum_{k=1}^n\int_0^{s_1}ds \int |\nabla Q_s g_{k,\lambda}|^2  d \mu_s .\label{Qs}
 \end{align}
 In what below we estimate the first and the second  terms in (\ref{Qs}),  respectively.

  For $s>0$,
 it follows from
  (\ref{gat}),  (\ref{boc}) and the right hand side of (\ref{Id}) that
\begin{align}
   \sum_{k=1}^n \int    |L_s( Q_sg_{k,\lambda} )|^2 d \mu_s
  \leq  &\sum_{k=1}^n \int   \Gamma_{2,s}( Q_sg_{k,\lambda},Q_sg_{k,\lambda} )  d \mu_s
  \label{abc}.
 \end{align}
Then, integration by parts, Cauchy-Schwarz inequality and (\ref{HEEE}) give
   \begin{align}
  &\Big|\sum_{k=1}^n \int \langle\nabla Q_s g_k , \nabla Q_sg_{k,\lambda}\rangle d \mu_s\Big| \nonumber\\
  =&\Big|\sum_{k=1}^n \int    (Q_s g_k)  L_s( Q_sg_{k,\lambda} ) d \mu_s\Big|
    &\nonumber\\
  \leq &\Big(\sum_{k=1}^n \int    |Q_s g_k |^2 d \mu_s\Big)^{1/2}\Big(\sum_{k=1}^n \int    |L_s( Q_sg_{k,\lambda} )|^2 d \mu_s\Big)^{1/2}\nonumber
  \\
\leq &  c_1^{-1/2} \big(s^{-1}\wedge 1)^{1/2} n^{1/2}\Big(\sum_{k=1}^n \int      \Gamma_{2,s}( Q_sg_{k,\lambda},Q_sg_{k,\lambda} )    d \mu_s\Big)^{1/2}
\nonumber
  \\
= & c_1^{-1/2} \big(s^{-1/2}\wedge 1)n^{1/2}  \Big( -\sum_{k=1}^n\frac{d}{ds}\int \Gamma_{1,s}( Q_sg_{k,\lambda},Q_sg_{k,\lambda} )d \mu_s \Big)^{1/2}   ,\label{gt}
 \end{align}
where   (\ref{di}) is used  with $k=1$ in the last step above.
Applying (\ref{gk}), we   have
 \begin{align}
  \sum_{k=1}^n\int\Gamma_{1,0}( g_{k,\lambda},g_{k,\lambda} )d\mu= & \sum_{k=1}^n\int |\nabla  g_{k,\lambda}|^2 d \mu\nonumber\\
= & -\int    (E_\lambda g_k) (L E_\lambda g_k ) d \mu\nonumber\\\leq &  \sum_{k=1}^n\lambda\int    | g_{k,\lambda}|^2 d \mu\nonumber\\
= &c_0^{-1}\lambda^2 n.\label{grad e}
 \end{align}
 Notice that $\log s_1+2\leq |\log\lambda|$ by $c_0<e^{-4}$.
 Then, it follows from (\ref{gt}),  (\ref{grad e}) and  Lemma \ref{Tecn}
    \begin{align}\label{grad epj}
 &\Big|\sum_{k=1}^n\int_0^{s_1} ds\int \langle\nabla Q_s g_k , \nabla Q_sg_{k,\lambda}\rangle d \mu_s\Big|\nonumber\\
 \leq &  c_1^{-1/2} n^{1/2}  \int_0^{s_1} \big(s^{-1/2 }\wedge 1)   \Big(-\sum_{k=1}^n\frac{d}{ds}\int \Gamma_{1,s}( Q_sg_{k,\lambda},Q_sg_{k,\lambda} )d \mu_s  \Big)^{1/2}   ds
  \nonumber\\
 \leq & c_0^{-1/2} c_1^{-1/2} \lambda | \ln \lambda |^{1/2} n.
 \end{align}

The gradient inequality   in  the proof of Lemma 2.5 of \cite{KP21} shows that  for $u \in H^1(\mu)$,
    \begin{align}
   \int |\nabla Q_s u |^2 d \mu_s\leq \int |\nabla  u|^2 d \mu ,\ \ \ \ \ s>0.\label{KP}
 \end{align}
When $u\in \mathcal{F}_\mu$, this inequality  is also a consequence of (\ref{di}) with $k=1$, (\ref{gat}) and the right hand side of
(\ref{Id}).
 We   have by (\ref{grad e}) and (\ref{KP})
   \begin{align}\label{grad epi}
 \int_0^{s_1}  \Big(\sum_{k=1}^n\int |\nabla Q_sg_{k,\lambda}|^2 d \mu_s\Big) ds \leq
&\int_0^{s_1}  \Big(\sum_{k=1}^n \int  |\nabla  g_{k,\lambda}|^2 d \mu \Big) ds
 \nonumber\\\leq &   c_0^{-1}\lambda^2n\int_0^{s_1}  ds\nonumber\\
 = & c_0^{-1/2} \lambda n  .
 \end{align}

Applying (\ref{Qs}),  (\ref{grad epj}), (\ref{grad epi}) and $\lambda\in (0,  e^{-1})$, we obtain
    \begin{align}
\Big| \sum_{k=1}^n\int \big(  Q_{s_1}(g_k-g_{k,\lambda}) \big) Q_{s_1}g_{k,\lambda}  d \mu_{s_1}\Big|
\leq &2c_0^{-1/2}c_1^{-1/2}   \lambda | \ln \lambda |^{1/2}n.\nonumber
   \nonumber
 \end{align}
This estimate and
  (\ref{HEEE}) show that
    \begin{align}
\sum_{k=1}^n\int   |Q_{s_1} g_{k,\lambda}|^2   d \mu_{s_1}\leq &\sum_{k=1}^n\int   |Q_{s_1} g_k |^2   d \mu_{s_1}+2\Big|
\sum_{k=1}^n\int \big(  Q_{s_1}(g_k-g_{k,\lambda}) \big)  Q_{s_1}g_{k,\lambda}  d \mu_{s_1}\Big| \nonumber\\
\leq & c_0^{-1/2}c_1^{-1}  \lambda  n+2c_0^{-1/2}c_1^{-1/2}   \lambda | \ln \lambda |^{1/2}n.
   \nonumber\\
\leq &3c_0^{-1/2}c_1^{-1/2}   \lambda | \ln \lambda |^{1/2}n.\nonumber
 \end{align}
Since  the   Rayleigh quotient of $g_{k,\lambda}$ at time zero is no more than $\lambda $ for $1\leq k\leq n$, we have by (\ref{control}), (\ref{gk}) and  $c_0<e^{-4}$
\begin{align}
 \sum_{k=1}^n \int |  Q_{s_1}g_{k,\lambda}|^2 d \mu_{s_1}\geq e^{- c_0^{1/2}}\sum_{k=1}^n\int |   g_{k,\lambda}|^2 d \mu\geq  (2c_0)^{-1}\lambda n .\nonumber
 \end{align}
Combing the two estimates above, we get
    \begin{align}
c_0^{-1}    \leq 6 c_0^{-1/2}c_1^{-1}| \ln \lambda |^{1/2},
   \nonumber
 \end{align}
   which implies the conclusion together with (\ref{gk}).
\qed\medskip

$\mathbf{Proof \ of\  Theorem\ \ref{theorem2}}$\

Let $\mu$ be an isotropic,     log-concave probability measure on $\mathbb{R}^n$.
Since  $\psi_n\geq \psi_1\geq  \sqrt{\pi/2}$, to prove
Theorem \ref{theorem2}, we only need to show that for some universal constant $C>0$,
\begin{align}\sigma_\mu \leq C |\log  \psi_\mu|+C,\ \ \ \ n\geq C.\label{re}\end{align}
Denote by $\lambda_\mu$ the
first nonzero eigenvalue of the self adjoint  operator $-L$  on $L^2(\mu)$. By the  Buser-Ledoux inequality from \cite{B82},\cite{L04},
$\psi_\mu^2\leq 9 \lambda_\mu^{-1}$. Therefore,
to prove (\ref{re}), it is enough to show that for some universal constant $C>0$,
\begin{align}\sigma_\mu \leq C |\log  \lambda_\mu|+C,\ \ \ \ n\geq C.\label{rer}\end{align}
Applying the approximation procedure in \cite{K23}, we can assume also that $\mu$ is regular.

Noticing that $a_0=0$, we have by  (\ref{one dir}) and (\ref{Step0})
  \begin{align}\label{medd}
   \frac{tn}{2} \leq \mathbb{E} |a_t|^2\leq  8tn,\ \ \ \ \ t\in [0,e^{-10^3}] ,\ \ n\geq C_3,
 \end{align}
where $C_3>0$ is a universal constant.
It follows from the fact that  $p_t(x)$ is a martingale and    Cauchy-Schwarz inequality that
  \begin{align}
 \mathbb{E}|a_t|^2\leq  \mathbb{E} \int |x|^2p_t(x)dx=\int |x|^2p_0(x)dx=n,\ \ \ \ t\geq 0,\ n\geq1.\nonumber
 \end{align}
Since $\mathbb{E} |a_t|^2$ is a increasing function of $t$ by (\ref{one dir}), the   estimate above, the  upper bound in  (\ref{medd})
and (\ref{Qq}) show that  (\ref{HEEE}) holds for $n\geq C_3$ and some universal constant $c_1\in(0,e^{-1})$. Then
Lemma \ref{TecN} yields
 \begin{align}
 F(\lambda)\leq  10^2 c_1^{-2}\lambda|\log\lambda|,\ \ \ \ \ \lambda\in (0,e^{-1}) ,\ n\geq C_3,\label{fl}
 \end{align}
 where
  \begin{align}
 F(\lambda)=n^{-1}\sum_{k=1}^n \int |E_\lambda x_k|^2 d\mu.\nonumber
 \end{align}

The following  formula is from  (48) in \cite{KL22} which is a direct consequence of the $H^{-1}$-inequality  (\ref{ct}):
 \begin{align}\label{thin}
 \sigma_\mu^2\leq &4\int_{\lambda_\mu}^\infty \lambda^{-2}F(\lambda)d\lambda.
 \end{align}
Notice  that $F(\lambda)\leq 1$ for $\lambda >0$. Then, for $n\geq C_3$,  we have by  (\ref{fl}) and    (\ref{thin})
  \begin{align}
 \sigma_\mu^2
   \leq &4\int_{e^{-1}}^\infty    \lambda^{-2}d\lambda  +4\cdot10^2 c_1^{-2} \int_{\lambda_\mu\wedge e^{-1}}^{e^{-1} }  \lambda^{-1} |\log \lambda| d\lambda
  \leq   2\cdot10^2 c_1^{-2}|\log (\lambda_\mu\wedge e^{-1})|^2,\nonumber
 \end{align}
 which  implies   (\ref{rer}).
 \qed\medskip

\begin{remark}\label{re2}
It is tempting to  erase the term $|\log\lambda|$ in the estimate (\ref{gk2})  which might  be possible.  Some straightforward calculations with the help of (\ref{Id}) and (\ref{boc}) may not be
enough   to give a well  control of the left hand side in (\ref{abc}).
\end{remark}

\

\

\

\noindent   Address:

\

\noindent Qingyang Guan

\

\noindent Institute of Applied Mathematics

\noindent Academy of Mathematics and Systems Science

\noindent Chinese Academy of Sciences

\noindent BeiJing 100190

\noindent  China

\

\noindent Email address: guanqy@amt.ac.cn


\begin{thebibliography}{123}


\bibitem{ABP03}  Anttila, M., Ball, K., Perissinaki, I.,  \textit{The central limit problem for convex
 bodies.} Trans. Amer. Math. Soc., Vol.  355, no. 12, (2003), 4723-4735.


\bibitem{BGL14} Bakry, D., Gentil, I., Ledoux, M., \textit{Analysis and Geometry of Markov Diffusion Operators.} Springer, Berlin,
2014.




\bibitem{BN12}  Ball, K., Nguyen, V. H., \textit{Entropy jumps for isotropic log-concave random vectors and spectral gap.} Studia
Math., Vol. 213, no. 1, (2012), 81-96.



\bibitem{BK20} Barthe, F., Klartag, B., \textit{Spectral gaps, symmetries and log-concave perturbations.} Bull. Hellenic Math.
Soc., Vol. 64, (2020), 1-31.


\bibitem{BK03}  Bobkov, S. G., Koldobsky, A., \textit{On the central limit property of convex
bodies.} Geometric aspects of functional analysis, Lecture Notes in
Math., 1807, Springer, Berlin, (2003), 44-52.



\bibitem{B86} Bourgain, J., \textit{On high-dimensional maximal functions associated to
convex bodies.} Amer. J. Math., Vol.  108, no. 6, (1986), 1467-1476.


\bibitem{B87} Bourgain, J., \textit{Geometry of Banach spaces and harmonic analysis.} Proceedings
 of the International Congress of Mathematicians, (Berkeley,
Calif., 1986), Amer. Math. Soc., Providence, RI, (1987), 871-878.



\bibitem{B91} Bourgain, J., \textit{On the distribution of polynomials on high dimensional convex sets.} Geom. Aspects of
Funct. Anal., Israel seminar (1989-90), Lecture Notes in Math., Vol. 1469, Springer, (1991), 127-137.







\bibitem{BL76} Brascamp, H.J., Lieb, E.H., 1976.  \textit{On extensions of the Brunn-Minkowski and
Pr$\acute{e}$kopa-Leindler theorems, including inequalities for log concave
functions, and with an application to the diffusion equation.}
  J. Funct. Anal. Vol.  22, no. 4, 366-389.

\bibitem{BGVV14} Brazitikos, S., Giannopoulos, A., Valettas, P., Vritsiou, B. -H.,
Geometry of isotropic convex bodies. American Mathematical Society, 2014.

\bibitem{B82}   Buser, P., \textit{A note on the isoperimetric constant.} Ann. Sci. $\acute{\mathrm{E}}$cole Norm.
Sup., Vol. 15, (1982), 213-230.

\bibitem{C21}   Chen, Y., \textit{An almost constant lower bound of the isoperimetric coefffcient in the KLS conjecture.} Geom.
Funct. Anal. (GAFA), Vol. 31, (2021), 34-61.





\bibitem{E13} Eldan, R., \textit{Thin shell implies spectral gap via a stochastic localization scheme.} Geom. Funct. Anal.
(GAFA), Vol. 23, (2013), 532-569.


\bibitem{EK11} Eldan, R.,  Klartag, B.,\textit{Approximately gaussian marginals and the hyperplane conjecture.}
Proc. of a workshop on "Concentration, Functional Inequalities and Isoperimetry", Contemporary
 Math., Vol. 545, Amer. Math. Soc., (2011), 55-68.





\bibitem{JLV22} Jambulapati, A., Lee, Y. T., Vempala, S., \textit{A slightly improved bound for the KLS constant.}, 2022.
arXiv:2208.11644.

\bibitem{KLS95} Kannan, R., Lov$\acute{\mathrm{a}}$sz, L., Simonovits M., \textit{Isoperimetric problems for convex bodies and a localization lemma.} Discr. Comput. Geom., Vol.13, (1995), 541-559.


\bibitem{K06} Klartag, B.,  \textit{On convex perturbations with a bounded isotropic constant.} Geom. and Funct. Anal. (GAFA), Vol.16, (2006), 1274-1290.


\bibitem{K09}
Klartag, B., \textit{A Berry-Esseen type inequality for convex bodies with an unconditional basis.} Probab. Theory
Related Fields, Vol. 145, (2009), 1-33.

\bibitem{K21} Klartag, B., \textit{On Yuansi Chen's work on the KLS conjecture.} Lecture notes prepared for a winter school
at the Hausdorff Institute, January 2021. Available at https://www.him.uni-bonn.de/fileadmin/
him/Lecture$\underline{\ }$Notes/chen$\underline{\ }$lecture$\underline{\ }$notes.pdf

\bibitem{K23}  Klartag, B., \textit{Logrithmic bounds for isoperimetry and slices of convex sets.} Ars Inveniendi Analytica (2023), Paper No.4, 17pp.



\bibitem{KP21} Klartag, B., Putterman, E., \textit{Spectral monotonicity under Gaussian convolution.} To appear in Ann. Fac.
Sci. Toulouse Math. arXiv:2107.09496.



\bibitem{KL22}  Klartag, B., Lehec, J., \textit{Bourgain's slicing problem and KLS isoperimetry up to polylog.} Geom. Funct. Anal.
(GAFA), Vol. 32, no. 5, (2022), 1134-1159.


\bibitem{KL24}  Klartag, B., Lehec, J., \textit{Isoperimetric inequalities in high-dimensional
convex sets.}  Lecture notes,
Institut Henri Poincar$\acute{\mathrm{e}}$ (IHP), Paris, (2024).



 \bibitem{KM22} Klartag, B., Milman, V., \textit{The slicing problem by Bourgain.} In: Analysis at Large, dedicated to the Life and
Work of Jean Bourgain, edited by A. Avila, M. Rassias and Y. Sinai, Springer, (2022), 203-232.







\bibitem{L04} Ledoux, M., \textit{Spectral gap, logarithmic Sobolev constant, and geometric
bounds.} Surveys in differential geometry.  Vol. IX, Int. Press, (2004), 219-240.


\bibitem{LV17} Lee, Y. T., Vempala, S., \textit{Eldan's stochastic localization and the KLS conjecture: Isoperimetry, Concentration
 and Mixing.} 58th Annual IEEE Symposium on Foundations of Computer Science (FOCS 2017),
IEEE Computer Soc., (2017), 998-1007.


\bibitem{LV24}  Lee, Y. T., Vempala, S.
\textit{Eldan's stochastic localization and the KLS conjecture: Isoperimetry, concentration and mixing.}  Ann. of Math., Vol. 199, no. 3, (2024),    1043-1092.










\end{thebibliography}
\end{document}